\documentclass[10pt,oneside,a4paper,reqno]{amsart}
\RequirePackage{fix-cm} 
\usepackage[T1]{fontenc}
\usepackage[english]{babel}
\usepackage{amssymb,amsmath,amsthm,graphicx,amsfonts}
\usepackage[foot]{amsaddr}
\usepackage{lmodern}

\usepackage[dvipsnames]{xcolor}
\definecolor{ourcolor}{RGB}{4,42,85}

\usepackage[
   rm={oldstyle=false,proportional=true},
   sf={oldstyle=true,proportional=true},
   tt={oldstyle=true,proportional=true,variable=true},
   qt=false,
 ]{cfr-lm}
\usepackage{cite,enumerate,setspace}
\usepackage[normalem]{ulem}
\usepackage{hyperref}
 \hypersetup{
 colorlinks,
 citecolor=ourcolor,
 linkcolor=ourcolor,
 urlcolor=black}

\numberwithin{equation}{section}

\usepackage{geometry}
 \makeatletter
 \def\@seccntformat#1{\hspace*{0mm}%
  \protect\textup{\protect\@secnumfont
    \ifnum\pdfstrcmp{subsection}{#1}=0 \bfseries\fi
    \csname the#1\endcsname
    \protect\@secnumpunct
      }%
 }

\newcommand{\assign}{:=}
\newcommand{\backassign}{=:}

\newcommand{\mathd}{\mathrm{d}}
\newcommand{\of}{:}

\newcommand{\tmabbr}[1]{#1}
\newcommand{\tmcolor}[2]{{\color{#1}{#2}}}
\newcommand{\tmem}[1]{{\em #1\/}}
\newcommand{\tmmathbf}[1]{\ensuremath{\boldsymbol{#1}}}

\newcommand{\tmname}[1]{\textsc{#1}}
\newcommand{\tmop}[1]{\ensuremath{\operatorname{#1}}}

\newcommand{\tmtextit}[1]{\text{{\itshape{#1}}}}

\newcommand{\tmtextup}[1]{\text{{\upshape{#1}}}}
\newenvironment{enumerateroman}{\begin{enumerate}[i.] }{\end{enumerate}}
\newtheorem{lemma}{Lemma}
\newtheorem{proposition}{Proposition}
{\theoremstyle{remark}\newtheorem{remark}{Remark}}
\newtheorem{theorem}{Theorem}

\newcommand{\Ddelta}{D_{\eps}}
\newcommand{\DMIdelta}{\mathcal{D}_{\eps}}
\newcommand{\Vdelta}{\mathcal{V}}
\newcommand{\Vtildedelta}{\tilde{\mathcal{V}}}

\newcommand{\myprop}[1]{#1}
\newcommand{\myspace}{H_{\eps} (\mathbb{R}^2;
\mathbb{S}^2)}
\newcommand{\myspacefull}{H_{\eps} (\mathbb{R}^2)}
\newcommand{\jmbar}{\overline{\tmmathbf{m}}}

\newcommand{\msq}{{\tmmathbf{m}_{\perp}^{\eps}}}
\newcommand{\bpi}{\tmmathbf{\pi}}
\newcommand{\nnn}{\tmmathbf{n}}
\newcommand{\uuu}{\tmmathbf{u}}
\newcommand{\eps}{\varepsilon}
\newcommand{\R}{\mathbb{R}}
\newcommand{\rr}{x}
\newcommand{\dd}{\mathd}
\newcommand{\sss}{y}

\newcommand{\m}{\tmmathbf{m}}
\newcommand{\mpa}{m_{\scriptscriptstyle{\|}}}
\newcommand{\mps}{m_{\scriptscriptstyle{\|}}^{\eps}}
\newcommand{\mpr}{\tmmathbf{m}_{\perp}}

\newcommand{\mpe}{\tmmathbf{m}_{\perp}}

\newcommand{\Om}{\Omega}

\newcommand{\Od}{\mathcal{O}^+_{\eps}}
\newcommand{\Odbar}{\mathcal{O}_{\bar{\eps}}}
\newcommand{\Odpm}{\mathcal{O}_{\eps}}

\newcommand{\Odpepsbar}{\mathcal{O}_{\bar{\eps}}}

\newcommand{\grad}{\nabla}
\newcommand{\divv}{\mathrm{div\,}}

\newcommand{\RR}{\mathbb{R}}

\newcommand{\Stwo}{\mathbb{S}}

\newcommand{\supp}{\textrm{\tmop{supp}}}

\newcommand{\jref}[1]{\tmtextup{\tmcolor{black}{\eqref{}}}}

\newcommand{\stt}{:}
\newcommand{\eqs}{=}

\newcommand{\dbar}{\tmcolor{black}{\bar{\eps}}}
\let\oldtocsection=\tocsection

\let\oldtocsubsection=\tocsubsection

\let\oldtocsubsubsection=\tocsubsubsection

\renewcommand{\tocsection}[2]{\hspace{0em}\oldtocsection{#1}{#2}}
\renewcommand{\tocsubsection}[2]{\hspace{1em}\oldtocsubsection{#1}{#2}}
\renewcommand{\tocsubsubsection}[2]{\hspace{2em}\oldtocsubsubsection{#1}{#2}}

\begin{document}

\title[{\sc Reduced energies for thin ferromagnetic films}]{{\large \sc \uppercase{Reduced energies for thin ferromagnetic films with perpendicular
anisotropy}}}
\author[1]{{\sc Giovanni Di Fratta}$^{(1)}$}
\address[1]{Dipartimento di Matematica e Applicazioni ``R.~Caccioppoli'', Università degli Studi di Napoli
    ``Federico II'', Via Cintia, 80126, Napoli, Italy.}
\author[2]{{\sc Cyrill B. Muratov}$^{(2,3)}$}
\author[3]{{\sc Valeriy V. Slastikov}$^{(4)}$}

\address[2]{Dipartimento di Matematica,
    Universit\`a di Pisa, Largo B. Pontecorvo, 5, 56127 Pisa,
    Italy.}   
\address[3]{Department of Mathematical Sciences, New Jersey
    Institute of Technology, Newark, NJ 07102, USA.}   
    
\address[4]{School of Mathematics, University of
Bristol, Bristol BS8 1TW, United Kingdom.
}


\begin{abstract}
  We derive four reduced two-dimensional models that describe, at
  different spatial scales, the micromagnetics of ultrathin
  ferromagnetic materials of finite spatial extent featuring
  perpendicular magnetic anisotropy and interfacial
  Dzyaloshinskii-Moriya interaction. Starting with a microscopic model
  that regularizes the stray field near the material's lateral edges,
  we carry out an asymptotic analysis of the energy by means of
  $\Gamma$-convergence. Depending on the scaling assumptions on the
  size of the material domain vs. the strength of dipolar interaction,
  we obtain a hierarchy of the limit energies that exhibit
  progressively stronger stray field effects of the material
  edges. These limit energies feature, respectively, a renormalization
  of the out-of-plane anisotropy, an additional local boundary penalty
  term forcing out-of-plane alignment of the magnetization at the
  edge, a pinned magnetization at the edge, and, finally, a pinned
  magnetization and an additional field-like term that blows up at the
  edge, as the sample's lateral size is increased. The pinning of the
  magnetization at the edge restores the topological protection and
  enables the existence of magnetic skyrmions in bounded samples.
 
  \vspace{4pt}
  
  \noindent {\scriptsize {\sc Mathematics Subject Classification.}
    82D40,49S05, 35C20.}
  
  \vspace{1pt}
  
  \noindent {\scriptsize {\sc Key words.} Micromagnetics, Magnetic
    thin films, Dzyaloshinskii--Moriya interaction,
    $\Gamma$-convergence.}
\end{abstract}

\begingroup
\def\uppercasenonmath#1{} 
\let\MakeUppercase\relax 
\maketitle
\endgroup

\tableofcontents

\section{Introduction and motivation}
With an ever-increasing control and sophistication of nanofabrication
techniques, there is a growing need for a better understanding of the
physical phenomena at the nanoscale that are determined by the material
geometry. In today's nano-magnetic systems, one typically encounters
materials consisting of one or several quasi-two-dimensional magnetic
layers interspersed with non-magnetic layers. The presence of
magnetic/non-magnetic interfaces gives rise to new physical effects
that have the potential to enable the next generation of
nanoelectronic devices that harness both the electric charge and the
spin degrees of freedom of electron for information technologies
\cite{soumyanarayanan16}. In the context of such {\em spintronics}
applications \cite{hirohata20,bhatti17,bader10}, one is particularly
interested in creating and manipulating spin configurations that are
endowed with non-trivial topological characteristics, which make them
robust against external influences and noise \cite{zhang20}.

The basic information unit in a topological spintronic device is the
{\em magnetic skyrmion} ---a particle-like continuous spin texture
with topological degree +1 (under a natural sign convention)
\cite{nagaosa13,braun12,fert13,fert17}. For this reason, there has
been considerable interest in the behavior of skyrmions in confining
geometries, both theoretically, computationally, and experimentally
(see, e.g., \cite{mmss:cmp23,rohart13,ohara21,boulle16,aranda18a,
  sampaio13,riveros21,cortes-ortuno19, Fratta2022, Davoli2022}; this
list is certainly not meant to be exhaustive).

However, under confinement, the topological protection of
non-collinear spin textures is a priori lost since the topological
degree of the spin configuration on a bounded domain is generally not
well-defined. In this case, a skyrmion-like spin texture may be
continuously deformed into a uniform magnetization state by pushing
the skyrmion out of the domain through the boundary.  A natural
solution is to settle the problem in the framework of curvilinear
magnetism. Indeed, magnetic thin films with the shape of closed
surfaces provide a concrete alternative for degree-preserving
confinements and, thus, toward the realization of chiral magnetic
textures.  The literature on this topic has grown very large. We refer
the reader to \cite{Fratta2023a,Fratta2023, DiFratta2020a,
  DiFratta2019, Gaididei2014, Sloika2017, Streubel2016, ignat21}, see also the
recent monograph \cite{Makarov2022}, for further reading on the
analysis of magnetic skyrmions in curved geometries that are close in
spirit to our interests here.  But as soon as one is interested in
planar thin films, further stabilization mechanisms for the magnetic
skyrmions in spintronic nanodevices would be required that provide a
repulsive interaction between the skyrmion and the device edge.

In this paper, we explore the additional energetic effects appearing
at the edges of two-dimensional ferromagnetic materials of finite
spatial extent. Due to the significant role played by the stray field
in ferromagnetic materials, these effects are often difficult to
predict and cause the emergence of new physical phenomena driven by
the material edges.  For example, in soft ferromagnets in the form of
thin films, the additional contribution of the stray field may
penalize the normal component of the magnetization at the film edge
\cite{CARBOU_2001, kohn05arma}, causing the appearance of boundary
vortices \cite{kurzke06,moser04,ignat23,official23}, edge-curling
domain walls \cite{lms:non18,lms:jns20}, interior walls
\cite{ignat10,ignat11a,mmns:arma23}, etc. (for a review, see
\cite{desimone06r}). In the current materials for spintronics
applications, additional physical mechanisms contribute at the film
edge \cite{rohart13,ms:prsla17,mskt:prb17}, further complicating the
situation. Therefore, to better understand the energetics of the
material edge, we carry out an asymptotic analysis of the
micromagnetic model of a two-dimensional ferromagnet exhibiting
perpendicular magnetic anisotropy and interfacial
Dzyaloshinskii-Moriya interaction \cite{rohart13}, using the
techniques of $\Gamma$-convergence. The singularity of the stray field
near the film edge in two-dimensional micromagnetics requires a
regularization of the magnetostatic problem near the edge and gives
rise to a hierarchy of reduced models appearing in the limit of the
vanishing stray field interaction strength. This regularization,
however, does not affect in any way the obtained limits, demonstrating
the universality of the asymptotic behavior of two-dimensional
ferromagnets.

We demonstrate that depending on the scaling of the lateral size of a
simply connected ferromagnetic domain with the strength of the
effective stray field interaction and for suitable renormalizations of
the other parameters, there are four distinct asymptotic regimes in
which the stray field acts differently at the domain edge and in the
ferromagnet's interior. In the first regime reminiscent of the thin
film limit studied by Gioia and James \cite{Gioia1997} (see also
\cite{CARBOU_2001, dmrs:sima20, DiFratta2020a}), the edge does not
exert any influence on the magnetization, resulting in a free boundary
condition and a renormalization of the magnetocrystalline anisotropy
constant (Theorem \ref{t:GJ}). In the second regime reminiscent of the
one studied by Kohn and Slastikov for soft ferromagnets
\cite{kohn05arma} and characterized by a larger lateral film extent,
the edge begins to exert an additional penalization of the deviation
of the magnetization from either one of the out-of-plane directions
(Theorem \ref{thm:thmKS}). In the third regime at yet larger film's
lateral extent, the magnetization becomes rigidly pinned to a single
out-of-plane direction at the film edge, while the stray field still
contributes locally in the film's interior (Theorem
\ref{thm:locengy}). Finally, in the fourth regime at yet larger film's
extent, the magnetization is also rigidly pinned at the film's edge,
but a nonlocal interaction term appears in the interior, as well as an
additional geometry-induced external field-like term (Theorem
\ref{thm:nonlocengy}). We note that the third regime corresponds
precisely to the one studied in \cite{mmss:cmp23}, where the topological
protection of single N\'eel skyrmions was shown to be restored in a
minimal micromagnetic setting.

The proofs of the main results are rather technical and require a
careful asymptotic analysis of the stray field close to the film
edge. For this purpose, we find it convenient to reformulate the
leading order expansion for the dipolar interaction energy
\cite{kmn:arma19} in terms of the gradients of the magnetization. This
is then used to compare the contributions of the mollified material
edge with that of the trace of the magnetization at a fixed material
boundary of the limit domain. The proofs proceed by a
divide-and-conquer strategy, in which the different parts of the
dipolar interaction are progressively isolated and ultimately
estimated by the limiting energy up to error terms that are
bounded by a small fraction of the exchange energy. The
$\Gamma$-limits are, in turn, proved by suitably combining the
different terms in the energy and passing to the limit directly. We
note that considerably finer estimates are required here to establish
our results compared to those needed in
\cite{Gioia1997, CARBOU_2001, kohn05arma, dmrs:sima20, DiFratta2020a}.

\subsection{Outline}
The paper is organized as follows. In section \ref{sec:model}, we 
introduce the micromagnetic model of a two-dimensional ferromagnetic
film exhibiting perpendicular magnetocrystalline anisotropy and an
interfacial DMI. Here we also introduce the edge regularization, which
is further derived from first principles in section
\ref{sec:micr-deriv-reduc}, and state our results informally, so that
they can be easily related to the original physical model. Then, in
section \ref{sec:discussion}, we state the precise assumptions and
definitions, and proceed to present the statements of the main results
of the paper. The rest of the paper is devoted to the proofs of the
theorems. In section \ref{sec:techlemmas}, we prove some preliminary
technical results. In section \ref{sec:asympt}, we carry out the
necessary asymptotic expansions of the magnetostatic energy. Finally,
in section \ref{sec:Gammaconv}, we complete the proofs of
$\Gamma$-convergence.

\subsection{Acknowledgments.} {\sc G.DiF.} acknowledges support from
the Austrian Science Fund (FWF) through the project {\emph{Analysis
    and Modeling of Magnetic Skyrmions}} (grant P34609). The work of
{\sc C.B.M.} was supported, in part, by NSF via grant DMS-1908709.
{\sc V.V.S.}  acknowledges support by Leverhulme grant
RPG-2018-438. {This study also received funding from the European
  Union – Next Generation EU – PRIN 2022 PNRR Project P2022WJW9H. {\sc
    C.B.M.}  is a member of INdAM/GNAMPA and acknowledges the MUR
  Excellence Department Project awarded to the Department of
  Mathematics, University of Pisa, CUP I57G22000700001.  {\sc G.DiF.}
  and {\sc C.B.M.}  thank the Hausdorff Research Institute for
  Mathematics in Bonn for its hospitality during the Trimester Program
  \emph{Mathematics for Complex Materials} funded by the Deutsche
  Forschungsgemeinschaft (DFG, German Research Foundation) under
  Germany Excellence Strategy – EXC-2047/1 – 390685813.  {\sc G.DiF.}
  also thanks TU Wien and MedUni Wien for their support and
  hospitality.  Finally, {\sc G.DiF.}  and {\sc V.V.S.} would like to
  thank the Max Planck Institute for Mathematics in the Sciences in
  Leipzig for support and hospitality.

\section{The micromagnetic model}\label{sec:model}
We start by considering a reduced two-dimensional micromagnetic energy
for an extended ferromagnetic thin film of effective thickness
$\delta > 0$, with the lengths measured in the units of the exchange
length $\ell = \sqrt{A / K_{\mathrm{d}}}$, where
$K_{\mathrm{d}} = \frac{1}{2} \mu_0 M_{\mathrm{s}}^2$, $\mu_0$ is the
vacuum permeability, $M_{\mathrm{s}}$ is the saturation magnetization
and $A$ is the exchange stiffness. The magnetization is characterized
by a non-dimensional vector $\m(x) \in \R^3$ at each point
$x \in \R^2$ of the film. We assume that the film exhibits
perpendicular magnetic anisotropy, interfacial Dzyaloshinskii-Moriya
interaction (DMI) and in the presence of an applied field
perpendicular to the film plane, so that the energy functional $E(\m)$
has the form
\begin{equation}
  E (\m) = E_{\mathrm{ex}} (\m) + E_{\mathrm{a}} (\m) + E_{\mathrm{Z}} (\m) +
  E_{\mathrm{DMI}} (\m) + E_{\mathrm{s}} (\m) . \label{eq:Einf}
\end{equation}
Here, in order of appearance, the terms are the exchange, anisotropy,
Zeeman, the DMI and the stray field energies measured in the units of
$A \ell \delta$.  As was discussed in
{\cite{ms:prsla17,bms:prb20,bms:arma21,kmn:arma19}}, in an extended
film where
\begin{equation}
  \m : \RR^2 \to \Stwo^2, \label{eq:mS2}
\end{equation}
is sufficiently smooth and goes to, say, $\m_0 = (0, 0, - 1)$ sufficiently
fast at infinity, with the notations
\begin{equation}
  \m = (\mpr, \mpa), \qquad \mpr : \R^2 \to \R^2, \quad \mpa : \R^2 \to \R,
  \label{eq:3}
\end{equation}
where $\mpr$ is the in-plane component and $\mpa$ is the out-of-plane
component of $\m$, respectively, these terms take the following form
\cite{bms:prb20,bms23}:
\begin{eqnarray}
  E_{\mathrm{ex}} (\m) & \assign & \int_{\R^2} | \nabla \m |^2  \dd \rr, \\
  E_{\mathrm{a}} (\m) & \assign & Q \int_{\R^2} | \mpr |^2  \dd \rr, \\
  E_{\mathrm{Z}} (\m) & \assign & - 2 h \int_{\R^2} (1 + \mpa)  \dd \rr, \\
  E_{\mathrm{DMI}} (\m) & \assign & \kappa \int_{\R^2} \left( \mpa  \divv 
  \mpr - \mpr \cdot \nabla \mpa \right)  \dd \rr, \\
  E_{\mathrm{s}} (\m) & \assign & - \int_{\R^2} | \mpr |^2  \dd \rr +
  \frac{\delta}{4 \pi} \int_{\RR^2} \int_{\R^2} \frac{\divv  \mpr \left( \rr
  \right) \cdot \divv  \mpr ( \sss)}{| \rr - \sss |}  \dd \rr  \dd
  \sss, \\
  &  & \qquad \qquad \qquad \qquad \qquad - \frac{\delta}{8 \pi} \int_{\R^2}
  \int_{\R^2} \frac{( \mpa ( \rr ) -  \mpa
  (\sss) )^2}{| \rr - \sss |^3}  \dd \rr  \dd \sss,  \label{eq:Es}
\end{eqnarray}
where
\begin{equation}
  Q = \frac{K_{\mathrm{u}}}{K_{\mathrm{d}}}, \qquad \kappa = \frac{D}{\sqrt{A
  K_{\mathrm{d}}}}, \qquad h = \frac{H}{M_{\mathrm{s}}}, \label{eq:params}
\end{equation}
with $Q$, $\kappa$ and $h$ being the dimensionless quality factor of
the out-of-plane anisotropy, the dimensionless DMI strength and the
dimensionless applied field strength, corresponding to the dimensional
magnetocrystalline anisotropy constant $K_{\mathrm{u}}$, DMI strength
$D$ normalized per unit volume, and the out-of-plane field $H$,
respectively. Note that $\kappa$ and $h$ may change sign, while for a
perpendicular magnetic anisotropy material we have $Q > 1$. Under
suitable conditions, the above energy exhibits local minimizers in the
form of the topologically non-trivial magnetization configurations --
\tmtextit{magnetic skyrmions}
{\cite{bogdanov89,bogdanov89a,bogdanov94,bogdanov94a,
    bogdanov99,ivanov90,melcher14,bms:prb20,bms:arma21,mmss:cmp23,bms23,
    Davoli2022}}.

Observe that the stray field energy in \eqref{eq:Es} admits the
following representation with the help of the Fourier transform
\begin{equation}
  \widehat{\m} (k) = \int_{\R^2} e^{- i k \cdot x}  (\m (x) - \m_0)  \dd \rr 
  \label{eq:mk}
\end{equation}
of $\m - \m_0 \in C^\infty_c (\R^2; \R^3)$:
\begin{equation}
  E_s (\m) = - \int_{\R^2} \left| \widehat{\m}_{\perp} (k) \right|^2 \frac{\dd
    k}{(2 \pi)^2} + \frac{\delta}{2} \int_{\R^2} \frac{|k \cdot
    \widehat{\m}_{\perp} (k) |^2}{| k |}  \frac{\dd k}{(2 \pi)^2} -
  \frac{\delta}{2} \int_{\R^2} | k |  | \widehat{m}_{\|} (k) |^2  \frac{\dd k}{(2
    \pi)^2} . \label{eq:Esk}
\end{equation}
In particular, the first term in the right-hand side of
\eqref{eq:Esk}, also referred to as the shape anisotropy term, may be
combined with $E_{\mathrm{a}} (\m)$ to define an effective
out-of-plane anisotropy with strength $Q - 1$ (going back to
\cite{winter61}); the second term in the right-hand-side of
\eqref{eq:Esk} represents the effect of the bulk charges and can be
seen to be non-negative; and the third term represents the effect of
the surface charges and is non-positive. The Fourier representation in
\eqref{eq:Esk} also arises as a relaxation of the energy in
\eqref{eq:Es} in the natural class of configurations in which
$\m - \m_0 \in H^1 (\R^2; \R^3)$. Notice that by \eqref{eq:Esk} and
simple interpolation inequalities the energy $E (\m)$ is always well
defined in this class (for further details, see {\cite{bms:arma21}}).

The expression for the stray field energy in \eqref{eq:Es} or
\eqref{eq:Esk} may be rigorously obtained as the leading order terms
in the asymptotic expansion of the full micromagnetic energy of a
three-dimensional ferromagnetic film of thickness $\delta \ll 1$, with
the errors controlled by the exchange energy at the next order
{\cite{kmn:arma19}}. It represents a suitably renormalized dipolar
interaction between different spins in an infinitesimally thin
ferromagnetic layer.  Care, however, is needed when extending this
definition to samples of \tmtextit{finite} spatial extent. Indeed, as
the stray field energy involves nonlocal terms, one cannot simply
restrict the integration in \eqref{eq:Es} to a {\tmem{bounded}}
spatial domain $\Om^\delta \subset \R^2$ (whose size may depend on
$\delta$), as this would disregard the dipolar interactions that
contribute to the first term in \eqref{eq:Es}. A more systematic
approach would instead consist of extending the magnetization
$\m : \Om^\delta \to \Stwo^2$ to the whole plane by zero outside
$\Om^\delta$. However, such an approach also presents difficulties, as
a jump discontinuity in $\m$ across $\partial \Om^\delta$ would then
generally make the nonlocal terms in \eqref{eq:Es} infinite.  Thus, a
regularization at the scale of the film thickness is necessary close
to the film edge to make sense of the energy in \eqref{eq:Es}. Such a
regularization was first introduced in {\cite{com:ieeetm13}} (see also
{\cite{mov:jap15}}, for further discussion see
{\cite{lms:non18,lms:jns20}}) in the context of reduced thin film
energies for soft ferromagnetic materials, in which the magnetization
tends to lie in the film plane.

We note that several regularizations are, in fact, possible that can
lead to slightly different reduced thin film energies. The precise
model would inevitably depend on the specific physics at the film
edge, which may be governed by a number of physical effects such as a
different material composition in an as-grown film near the edge,
changes in the crystalline structure near the edge, edge roughness,
etc. We point out, however, that the magnetization, which in the
physical space rotates on the scale of the exchange length that
exceeds by an order of magnitude the atomic scale {\cite{hubert}},
should experience the effect of the edge via some sort of effective
boundary terms. This is indeed confirmed by rigorous studies of the
thin film limit of soft-three dimensional ferromagnetic layers
{\cite{kohn05arma}}. The present paper aims to derive these
boundary terms via $\Gamma$-convergence for ferromagnetic films with
perpendicular magnetocrystalline anisotropy that are relevant to the
studies of magnetic skyrmions.

Our starting point will be the regularization in which for, say, a
given $\m \in C^\infty(\R^2; \Stwo^2)$ we define the physically
observable magnetization $\m_{\delta}$ in the film:
\begin{equation}
  \m_{\delta} ( \rr ) \assign \eta_{\delta} ( \rr )  \m
  ( \rr ) \qquad \forall \rr \in \R^2, \label{eq:meta}
\end{equation}
where for a bounded open, simply connected set $\Om^\delta$ with
boundary of class $C^2$ representing the film of finite extent we
defined a cutoff function
\begin{equation}
  \eta_{\delta} (\rr) \assign \eta \left( \frac{d_{\partial \Om^\delta}
      (\rr)}{\delta} \right), \label{eq:etadelta}
\end{equation}
in which $d_{\partial \Om^\delta}$ is the signed distance from the
boundary ({\tmabbr{cf.}}~\eqref{eq:signdist}) and $\eta$ is a
non-increasing, sufficiently regular function that goes from
$\eta(-\infty) = 1$ to $\eta(+\infty) = 0$ (see section
\ref{sec:discussion} for details). For a microscopic derivation of
this condition, see section~\ref{sec:micr-deriv-reduc}. Notice that
$\m_{\delta}$ thus defined automatically lies in $H^1 (\R^2; \R^3)$ if
$\m \in H^1_\mathrm{loc} (\R^2; \Stwo^2)$. We then replace all the
instances of $\m$ in \eqref{eq:Einf} with $\m_{\delta}$ to define the
reduced thin film energy
$\mathsf E_{\delta} (\m) \assign E (\m_{\delta})$.

We next specify the asymptotic regimes in which the obtained energy
$\mathsf E_{\delta}$ can be significantly simplified. In the first two
regimes the limit energy becomes local, with the edge effects either
disappearing or appearing as a boundary term that generalizes the
regime for soft ferromagnetic films identified by Kohn and Slastikov
{\cite{kohn05arma}}. We will also identify the scalings of the
parameters for which the resulting limit energy still exhibits the
terms that are needed to produce skyrmion-type solutions. To this end,
we introduce a small parameter $\eps > 0$ and make all the model
parameters, as well as the domain, depend on $\eps$ as follows:
\begin{equation}
  Q_{\eps} = 1 + \frac{\eps | \ln \eps |}{2 \pi \gamma_\eps}
  \hspace{0.17em} \alpha, \quad h_{\eps} = \frac{\eps | \ln \eps |}{2
    \pi \gamma_\eps} \hspace{0.17em} \beta, \quad \kappa_{\eps} = \left(
    \frac{\eps | \ln \eps |}{2 \pi \gamma_\eps} \right)^{1 / 2} \lambda,
  \quad \delta_{\eps} = \left( \frac{2 \pi \eps \gamma_\eps}{| \ln \eps |}
  \right)^{1 / 2}, \label{eq:parameps} 
\end{equation}
together with
\begin{equation}
  \Om^\delta_\eps \assign \eps^{- 1} \delta_{\eps}  \Om, \label{eq:Omeps}
\end{equation}
for some fixed $\lambda > 0$, $\alpha, \beta \in \R$ and
$\Om \subset \R^2$. Notice that at this point \eqref{eq:parameps} and
\eqref{eq:Omeps} simply represent a reparametrization that imposes a
certain dependence on two parameters, $\eps$ and $\gamma_\eps$ (with
the dependence of the latter on $\eps$ to be specified), and forces
suitable balances between different terms in the energy as
$\eps \to 0$ depending on the choices of $\gamma_\eps$.

With the above choices and after some algebra, we have
$\mathsf E_{\delta_\eps} (\m (\eps^{- 1} \delta_{\eps} \cdot)) =
E_{\eps} (\m) + C_{\eps}$, where
\begin{eqnarray}
  E_\eps (\m)
  & \assign & \int_{\R^2} \left( \eta_{\eps}^2 | \nabla \m |^2 +
              \alpha \eta_{\eps}^2 | \mpr |^2 - 2 \beta \eta_{\eps}
              \mpa \right)  \dd \rr 
              \nonumber\\ 
  &  & \quad + \lambda \int_{\R^2} \eta_{\eps}^2 \left( \mpa  \divv
       \mpr - 
       \mpr \cdot \nabla \mpa \right)  \hspace{0.17em} \dd \rr
       \nonumber\\ 
  &  & \quad + \frac{\gamma_\eps}{2 | \ln \eps |} \int_{\R^2}
       \int_{\R^2} 
       \frac{\divv (\eta_{\eps} \mpr) ( \rr )  \divv  (\eta_{\eps} \mpr)
       (\sss)}{| \rr - \sss |}  \dd \rr  \dd \sss \nonumber\\
  &  & \quad - \frac{\gamma_\eps}{4 | \ln \eps |} \int_{\R^2}
       \int_{\R^2} 
       \frac{(\eta_{\eps} ( \rr ) \mpa ( \rr ) - \eta_{\eps}
       ( \sss) \mpa ( \sss))^2}{| \rr - \sss |^3}  \dd \rr 
       \dd \sss,  \label{defmme} \label{eq:Eeps}
\end{eqnarray}
and the additive constant $C_{\eps}$ is independent of $\m$ and,
therefore, is inconsequential for the variational problem associated
with $\mathsf E_{\delta_\eps}$. We note that for $\lambda = 0$, only a
slightly different version of this type of energy with $\eps \sim 1$
can be shown to arise from the full micromagnetic energy of a
three-dimensional thin ferromagnetic film with variable thickness
equal to $\eta_{\eps} ( \rr )$, which tapers off at the film edge
{\cite{slastikov05}}.

Now, with the choice $\gamma_\eps \to \gamma$ as $\eps \to 0$ for
$\gamma > 0$, the limit functional will be shown to be
\begin{align}
  F (\m)
  & \assign \int_{\Om} \left( \left| \nabla \m \right|^2 + \alpha |
    \mpr |^2 - 2 \beta \mpa \right)  \dd \rr \nonumber\\
  &   + \lambda \int_{\Om} \left( \mpa  \divv  \mpr -
    \mpr \cdot 
    \nabla \mpa \right)  \dd \rr + \gamma \int_{\partial
    \Om} \left( (\mpr \cdot 
    \tmmathbf{n})^2 - \mpa^2 \right)  \dd \mathcal{H}^1
    (\sigma),  \label{eq:F} 
\end{align}
where $\tmmathbf{n}$ is the outward unit normal to $\partial \Om$, and
$F(\m)$ is defined for $\m \in H^1 (\Om; \Stwo^2)$. The limit
functional sees only the limit domain $\Omega$ and in addition to the
expected local terms inside $\Omega$ it features a boundary term that
penalizes the deviations of the in-plane component of the
magnetization from tangential to $\partial \Omega$, and another
boundary term that favors the out-of-plane component of the
magnetization to be $\pm 1$. These terms arise, respectively, as the
limits of the next-to-last and the last term in the definition of
$E_\eps$ in \eqref{eq:Eeps} due to the logarithmic divergence of the
respective integrals as $\eps \to 0$. We remark that the limit energy
with $\gamma = 0$ similarly arises when $\gamma_\eps \to 0$ and
$\eps \to 0$, a result analogous to a well-known result of Gioia and
James \cite{Gioia1997}.

We will also consider two other scaling regimes, which lead to
different limit behaviors. First, we define
\begin{equation}
  E_{\eps}^0 (\m) \assign E_{\eps} (\m) + \gamma_{\eps} \mathcal{H}^1
  (\partial \Om), \label{eq:Eepsinf}
\end{equation}
where $\gamma$ is replaced with $\gamma_{\eps}$ in \eqref{eq:Eeps},
and we will be interested in the limit in which
$\gamma_{\eps} \to + \infty$ as $\eps \to 0$ with $\alpha, \beta$ and
$\lambda$, as well as the domain $\Om$, fixed (note that the limit
$\gamma_\eps \to 0$ is much simpler and is obtained by just setting
$\gamma = 0$ in \eqref{eq:F}). When $\gamma_\eps \to +\infty$, we show
that for $\gamma_{\eps} = o (| \ln \eps |)$ the limit energy for
$E_{\eps}^0$ is given by
\begin{equation}
  F_0 (\m) \assign \int_{\Om} \left( \left| \nabla \m \right|^2 + \alpha
  \left| \mpr \right|^2 - 2 \beta \mpa \right)  \dd \rr + \lambda \int_{\Om}
  \left( \mpa  \divv  \mpr - \mpr \cdot \nabla \mpa \right)  \dd \rr,
  \label{eq:F0}
\end{equation}
specified for all $\m \in H^1 (\Om; \Stwo^2)$ such that
$\m =\tmmathbf{e}_3$ or $\m = -\tmmathbf{e}_3$ on $\partial \Om$ in
the sense of trace. This could be thought of in some sense as the
limit case of the energy in \eqref{eq:F} with $\gamma = \infty$, after
a suitable renormalization.

Finally, we consider the regime in which for $\nu > 0$ and
$\alpha, \beta, \lambda$ real we have
\begin{equation}
  Q_{\eps} = 1 + \frac{\eps}{2 \pi \nu} \alpha, \quad h_{\eps} =
  \frac{\eps}{2 \pi \nu}  \hspace{0.17em} \beta, \quad \kappa_{\eps} =
  \left( \frac{\eps}{2 \pi \nu} \right)^{1 / 2} \lambda, \quad
  \delta_{\eps} = \left( 2 \pi \eps \nu \right)^{1 / 2},
  \label{eq:paramepsnu}
\end{equation}
once again together with \eqref{eq:Omeps}, which corresponds to the
choice of $\gamma_{\eps} = \nu | \ln \eps |$ in
\eqref{eq:Eepsinf}. Here we find the following limit energy for
$E_\eps^0$, up to an additive constant:
\begin{align}
  F_{\nu} (\m)
  & \assign \int_{\Om} \left( \left| \nabla \m \right|^2 +
    \alpha \left| \mpr \right|^2 - 2 \beta  \mpa \right)  \dd \rr
    \notag \\ 
  & + \lambda
    \int_{\Om} \left( \mpa  \divv  \mpr - \mpr \cdot \nabla \mpa
    \right) \dd \rr + \nu \int_\Om \tmmathbf{b} \cdot \nabla
    \mpa   \label{eq:Fnu} 
    \dd \rr
  \\
  &   + \frac{\nu}{2}  \int_{\Om} \int_{\Om} \frac{\divv  \mpr
    ( \rr )  \divv  \mpr (\sss)}{| \rr - \sss |}  
    \dd \rr  \dd \sss - \frac{\nu}{2}  \int_{\Om} \int_{\Om}
    \frac{\nabla \mpa
    ( \rr ) \cdot \nabla \mpa (\sss) }{| \rr - \sss |}  \dd \rr  \dd
    \sss,  \notag
\end{align}
where we defined the vector field
\begin{align}
  \tmmathbf{b}(x) :=\int_{\partial \Omega}
  \frac{\mpa(y) \tmmathbf{n} (y)}{|x - y|} \dd y \qquad x \in
  \Om,   \label{eq:b} 
\end{align}
in which $\tmmathbf{n}$ is the outward unit normal to $\partial
\Om$. This vector field encodes the stray field effect of the film
edge. The energy $F_\nu$ is defined for all
$\m \in H^1 (\Om; \Stwo^2)$ such that $\m =\tmmathbf{e}_3$ or
$\m = -\tmmathbf{e}_3$ on $\partial \Om$ in the sense of trace. Note
that for $\m = \pm \tmmathbf{e}_3$ on $\partial \Omega$ we have
$\tmmathbf{b} \in C^\infty(\Omega)$ and $\tmmathbf{b}(x)$ diverges
logarithmically with distance as $x \in \Omega$ approaches
$\partial \Omega$. In particular, the term in the energy involving
$\tmmathbf{b}$ is under control by the gradient squared term.

All of the aforementioned statements are made precise within the framework of
$\Gamma$-convergence in section \ref{sec:discussion}.

\subsection{A microscopic derivation of the reduced two-dimensional
  model}\label{sec:micr-deriv-reduc}

As was already mentioned, the precise behavior of the magnetization
near the film edge depends on the detailed physics at the edge of the
film. Here we use a particular model that illustrates how an energy of
the form given in \eqref{eq:Eeps} may be obtained from a more
microscopic description.

To avoid dealing with truly discrete models of ferromagnetism at the
atomic scale, we pick a model that still allows to describe the film
as a continuum, but retains the thermodynamic essence of the
ferromagnetic phase and allows to evaluate the additional effects of
the film edge. Namely, we consider a mean-field model of a Heisenberg
ferromagnet with a long-range Kac attractive interaction. Such models
have been rigorously derived in the context of theories of phase
transitions, going back to Lebowitz and Penrose for the liquid-gas
phase transition {\cite{lebowitz66}} and Thompson and Silver for the
classical Heisenberg magnet {\cite{thompson73}}. Moreover, in the
considered limit the metastable spatially varying states may be
understood via minimization of a free energy functional
{\cite{giacomin96}}, which in the case of the Heisenberg model with
the interaction kernel $J_{\delta} (\left| \rr \right|)$ takes the
form
\begin{align}
  \mathcal F (\rho)
  & \eqs  - \frac{1}{2}  \int_{\Stwo^2} \int_{\Stwo^2} \int_{\Om}
           \int_{\Om} J_{\delta} \left( \left| \rr - \sss \right| \right) \left( \m
           \cdot \m' \right) \rho \left( \rr, \m \right) \hspace{0.17em} \rho \left( y,
           \m' \right)  \dd \rr  \, \dd \sss  \hspace{0.17em} \dd \mathcal{H}^2 \left( \m
           \right)  \dd \mathcal{H}^2 (\m') \nonumber\\
  &   \qquad  \qquad \qquad \qquad + \beta^{- 1} 
       \int_{\Stwo^2} \int_{\Om} \rho (\rr, \m) \ln \rho (\rr, \m) \hspace{0.17em}
       \dd \rr  \, \dd \mathcal{H}^2 \left( \m \right)
       .  \label{eq:FHeis} 
\end{align}
Here $\rho \in L^1 (\R^2 \times \Stwo^2 ; [0, \infty])$ is the probability
density to observe a spin at point $\rr \in \R^2$ in the direction $\m \in
\Stwo^2$, $J_{\delta} \in C^{\infty}_c (\R)$ is a positive, even interaction
potential such that
\begin{equation}
  \mathrm{supp} ( J_{\delta}) \subset B_{\delta} (0) \quad \text{and} \quad
  \int_0^{\infty} 2 \pi r J_{\delta} (r)  \dd r = J_0 > 0
\end{equation}
for a $J_0 > 0$ fixed independently of $\delta$, and $\beta > 0$ is the
inverse temperature. The function $\rho$ satisfies the following normalization
conditions:
\begin{equation}
  \int_{\Stwo^2} \rho (\rr, \m) \hspace{0.17em} \dd \mathcal{H}^2 (\m) = 1
  \quad \text{if} \; \rr \in \Om, \quad \rho \left( \rr, \m \right) = 0 \quad
  \text{if} \; \rr \in \RR^2 \setminus \Om, \label{eq:rhonorm}
\end{equation}
expressing the fact that the ferromagnet occupies the spatial domain $\Om
\subset \R^2$. For our purposes, all other terms in the energy, which are all
small perturbations to the Heisenberg exchange, have been neglected. Notice
that the parameter $\delta$ measures the finite range of the ferromagnetic
coupling and physically corresponds to the extent of the exchange interaction
of several lattice spacings.

The free energy in \eqref{eq:FHeis} admits a moments closure, allowing to
reduce the minimization problem to that of a functional of the average
magnetization (see also {\cite{fatkullin05}})
\begin{equation}
  \jmbar ( \rr ) \assign \int_{\Stwo^2} \m \rho (\rr, \m) 
  \hspace{0.17em} \dd \mathcal{H}^2 (\m) . \label{eq:mbar}
\end{equation}
For a fixed value of $\jmbar$ the entropic term in the free energy is
easily seen to be minimized pointwise by
\begin{equation}
  \bar{\rho} (\rr, \m) = \exp \left( \beta \left( \mu ( \rr )
  +\tmmathbf{\lambda} ( \rr ) \cdot \m \right) \right),
  \label{eq:rhob}
\end{equation}
where the functions $\mu$ and $\tmmathbf{\lambda}$ are obtained by enforcing
\eqref{eq:rhonorm} and \eqref{eq:mbar} with $\rho = \bar{\rho}$ in $\Om$:
\begin{eqnarray}
  1 & \eqs & 4 \pi e^{\beta \mu} \; \frac{\sinh (\beta
             |\tmmathbf{\lambda}|)}{\beta |\tmmathbf{\lambda}|}, \\
  \jmbar & \eqs & 4 \pi e^{\beta \mu} \; \frac{\beta |\tmmathbf{\lambda}|
                  \cosh (\beta |\tmmathbf{\lambda}|) - \sinh (\beta
                  |\tmmathbf{\lambda}|)}{\beta^2
                  |\tmmathbf{\lambda}|^3} \tmmathbf{\lambda}.  
\end{eqnarray}
This yields $\tmmathbf{\lambda}= \beta^{- 1}  \jmbar f (| \jmbar |) / | \jmbar
|$, where the function $f (s) \geqslant 0$ is the unique positive solution of
the equation
\begin{equation}
  s = \tmop{coth} f(s) - \frac{1}{f (s)}, \qquad 0 < s < 1, \label{eq:f}
\end{equation}
vanishing at $s = 0$ and diverging as $s \to 1^-$. The plot of $f (s)$
is presented in Fig.~\ref{fig:fU}(a). Note that
$f \in C^{\infty} ([0, 1))$ and is strictly monotone
increasing. Substituting this back to the entropy term results in
\begin{equation}
  \int_{\Stwo^2} \bar{\rho} (\rr, \m) \ln \bar{\rho} (\rr, \m)
  \hspace{0.17em} \dd \mathcal{H}^2 (\m) \eqs \ln \left( \frac{f (| \jmbar
  ( \rr ) |)}{4 \pi \sinh f (| \overline{\m} ( \rr ) |)}
  \right) + \jmbar ( \rr ) f (\left| \jmbar ( \rr )
  \right|) . \label{eq:entr}
\end{equation}
\begin{figure}[t]
  \begin{center}
    \resizebox{15cm}{!}{\includegraphics{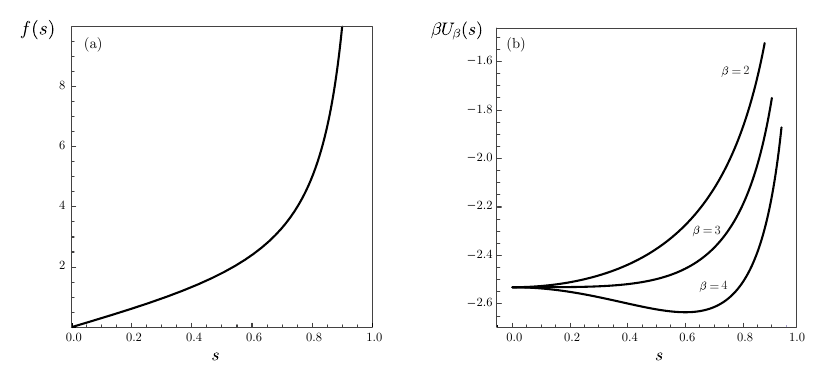}}
  \end{center}
  \caption{(a) Plot of $f (s)$. (b) Plots of $\beta U_{\beta} (s)$ for the
  indicated values of $\beta$ and $J_0 = 1$.\label{fig:fU}
  }
\end{figure}{\noindent}Thus, for $\jmbar ( \rr )$ fixed the free
energy satisfies
$\mathcal F (\rho) \geqslant \bar{\mathcal F} (\jmbar)$, where
\begin{align} \bar{\mathcal F} (\jmbar) \assign - \frac{1}{2} \int_{\R^2}
  \int_{\R^2} J_{\delta} \left( \left| \rr - \sss \right| \right)
  \jmbar ( \rr ) \cdot \jmbar ( \sss) \dd \rr \dd \sss + \frac{J_0}{2}
  \int_{\R^2} | \jmbar |^2 \dd \rr + \int_{\R^2} U_{\beta} (| \jmbar
  |) \dd \rr,
   \label{eq:Fbar} \end{align}
with equality holding if and only if $\rho (\rr, \m) = \bar{\rho} (\rr, \m)$.
Here the effective potential $U_{\beta}$ is given by
\begin{equation}
  U_{\beta} (s) \assign \beta^{- 1} \ln \left( \frac{f (s)}{4 \pi \sinh \, f
  (s)} \right) + \beta^{- 1} s f (s) - \frac{1}{2} J_0 \, s^2, \label{eq:U}
\end{equation}
with the convention that $U_{\beta} (0) \assign - \beta^{- 1} \ln 4 \pi$. The
plots of $U_{\beta} (s)$ for several values of $\beta$ are illustrated in
Fig.~\ref{fig:fU}(b).

Notice that the reduced energy in \eqref{eq:Fbar} may be rewritten in a more
convenient form as
\begin{equation}
  \bar{\mathcal F} (\jmbar) = \frac{1}{4}  \int_{\R^2} \int_{\R^2} J_{\delta} \left(
  \left| \rr - \sss \right| \right)  \left| \jmbar ( \rr ) - \jmbar
  ( \sss) \right|^2  \dd \rr  \dd \sss + \int_{\R^2} U_{\beta} (|
  \jmbar |)  \dd \rr, \label{eq:Fbar2}
\end{equation}
which is a vectorial, nonlocal analog of the classical Cahn-Hilliard
functional, since $U_{\beta}$ has a form of a Mexican hat potential
for $\beta > \beta_c \assign 3 J_0$. It is also easy to see that in a
periodic setting the energy functional $\bar{\mathcal F}$ admits a
unique minimizer $\jmbar =\tmmathbf{0}$ whenever
$\beta \leqslant \beta_c$, and a family of minimizers
$| \jmbar | = s_0 (\beta)$ with $0 < s_0 (\beta) < 1$ for
$\beta > \beta_c$ (see also {\cite{fatkullin05}}). To simplify matters
further, we employ the usual gradient approximation to the nonlocal
term in \eqref{eq:Fbar2} to obtain
$\bar{\mathcal F} (\jmbar) \simeq \bar{\mathcal F}_0 (\jmbar)$, where
(see also {\cite{fatkullin08}} for a closely related problem)
\begin{equation}
  \bar{\mathcal F}_0 (\jmbar) \assign \int_{\Om} \left( g_{\delta} |
    \nabla \jmbar |^2 + 
    U_{\beta} (| \jmbar |)  \right) \dd \rr, \label{eq:Fbar0}
\end{equation}
and $g_{\delta} \assign \frac{\pi}{4} \int_0^{\infty} r^3 J_{\delta} (r)  \dd
r = O (\delta^2)$. The expression in \eqref{eq:Fbar0} is specified for $\jmbar
\in H^1_0 (\Om; \R^3)$, inheriting the zero boundary condition from the
assumption that $\jmbar =\tmmathbf{0}$ in $\R^2 \setminus \Om$.

Near the edge of the sample the curvature of the edge is negligible to the
leading order in $\delta$. Hence, the problem of minimizing the energy in
\eqref{eq:Fbar0} reduces to a one-dimensional problem on half-line, i.e., to
minimizing the energy
\begin{equation}
  \bar{\mathcal F}_0^{\mathrm{1d}} (\jmbar) := \int_0^{\infty} \left( g_{\delta} |
    \jmbar' |^2 + U_{\beta} (| \jmbar |) - U_{\beta} (s_0 (\beta)) \right)  \dd
  \rr \label{eq:F01d}
\end{equation}
over
$\jmbar \in H^1_{\tmop{loc}} (\R^+; \R^3) \cap C (\overline{\R}^+;
\R^3)$ such that $\jmbar (0) =\tmmathbf{0}$. An explicit
energy-minimizing profile may be obtained from \eqref{eq:F01d}, using
the polar representation $\jmbar (x) = \phi (x) \tmmathbf u (x)$,
where $| \tmmathbf u | = 1$, for which we get
\begin{equation}
  \bar{\mathcal F}_0^\mathrm{1d} (\jmbar) = \int_0^{\infty} \left( g_{\delta} | \phi' |^2 +
  U_{\beta} (\phi) - U_{\beta} (s_0 (\beta))  \right) \dd x + \int_0^{\infty}
  g_{\delta} \phi^2 | \tmmathbf u ' |^2  \dd x.
\end{equation}
Thus, the energy $\bar{\mathcal F}_0^{\mathrm{1d}}$ is minimized by
$\tmmathbf u = \mathrm{const}$ and $\phi = \phi_\delta$, where by the
Modica-Mortola trick {\cite{modica87}} we have
\begin{equation}
  \min \bar{\mathcal F}_0^{\mathrm{1d}} = \bar{\mathcal
    F}_0^{\mathrm{1d}}  (\phi_{\delta} \tmmathbf u) = 2 \int_0^{s_0
    (\beta)} \sqrt{g_{\delta} (U_{\beta} (\phi) - 
    U_{\beta} (s_0 (\beta)))}  \hspace{0.17em} \dd \phi,
\end{equation}
\begin{equation}
  \int_0^{\phi_{\delta} (x)} \frac{\dd \phi}{\sqrt{g_{\delta}  (U_{\beta}
  (\phi) - U_{\beta} (s_0 (\beta)))}} = x, \qquad \tmmathbf u \in \Stwo^2 \
  \text{arbitrary} .
\end{equation}
In particular, we have $| \jmbar (x) | = \phi_{\delta} (x)$, where
$\phi_\delta(x)$ is a monotone increasing function such that
$\phi_\delta(0) = 0$ and $\phi'_\delta(0) > 0$, and $\phi_\delta(x)$
approaches exponentially the ``saturation magnetization''
$s_0 (\beta)$ for $x \gg \delta$. Therefore, as a matter of modeling
convenience one could replace the function $\phi_\delta$ with that of
$s_0 \eta_\delta$, where $\eta_\delta$ is defined in the following
section.

\section{Mathematical setup and statement of the main
  results}\label{sec:discussion}

\noindent We now give the precise mathematical formulation of the
considered problem and state our main theorems.

\subsection{Film geometry} \label{sec:geom}

Throughout the paper we assume that $\Om \subset \RR^2$ is a bounded,
simply connected open set with boundary of class $C^2$. We parametrize
$\partial \Om$ by a $C^2$ regular curve, and we denote by
\begin{align} \varphi : s \in I_{\partial \Om} \mapsto \varphi (s) \in \partial
  \Om, \quad I_{\partial \Om} \assign \left[ 0, \mathcal{H}^1 \left(
      \partial \Om \right) \right], \end{align} its positive parameterization
by the arc length (extended periodically, if necessary). We denote by
$\left( \tmmathbf{\tau} ( \rr ), \nnn ( \rr ) \right) \in \Stwo^1
\times \Stwo^1$ the respective Frenet frame at $\rr \in \partial
\Om$. The vector fields
$\tmmathbf{\tau}: \partial \Om \rightarrow \Stwo^1$ and
$\nnn : \partial \Om \rightarrow \Stwo^1$ represent, respectively, the
unit tangent vector field to $\partial \Om$ given by
$\tmmathbf{\tau} (\varphi(s)) \assign \varphi' (s)$, and the outer
unit normal vector field. To avoid cumbersome notations, we write
$\tmmathbf{\tau} (s)$ and $\nnn (s)$ to mean the compositions
$\tmmathbf{\tau} (\varphi (s))$ and $\nnn (\varphi (s))$ from now
on. Clearly, for every $s \in I_{\partial \Om}$ we have
$\tmmathbf{\tau} (s) \cdot \nnn (s) = 0$, and the Frenet-Serret
formulas hold:
\begin{align}
  \tmmathbf{\tau}' (s) & \eqs  - \kappa (s) \nnn (s),  \label{eq:Frenet1}\\
  \nnn' (s) & \eqs  \kappa (s) \tmmathbf{\tau} (s),  \label{eq:Frenet2}
\end{align}
where $\kappa (s) \assign - \varphi'' (s) \cdot \nnn (s)$ stands for
the signed curvature of $\partial \Om$ at the point $\varphi (s)$. Note
that since we assume $\partial \Om$ of class $C^2$, the vector fields
$\tmmathbf{\tau} (s)$ and $\tmmathbf{n} (s)$ are of class
$C^1 \left( I_{\partial \Om}; \Stwo^1 \right)$ and $\kappa (s)$ is a
continuous function.

In order to define a cutoff function near the boundary of the
ferromagnet, we first introduce the signed distance function from
$\partial \Om$ which assigns positive values to points in the exterior
of $\Om$ and negative values in the interior of $\Om$:
\begin{equation}
  d_{\partial \Om} (\rr) \assign \left\{ \begin{array}{ll}
                                           + \inf_{\sss \in \partial
                                           \Om} | \rr - \sss |
                                           & \text{if } \rr \in \RR^2
                                             \setminus \Om,\\
                                           - \inf_{\sss \in \partial
                                           \Om} | \rr - \sss |
                                           & \text{if } \rr \in \Om .
  \end{array} \right. \label{eq:signdist}
\end{equation}
Since $\Om$ is a $C^2$-domain, there exists $\dbar > 0$ such that for
any $0 < \eps < \dbar$ the set
\begin{align}
  \Odpm \assign \left\{ \rr \in \RR^2 \hspace{0.17em} \stt \left|
  d_{\partial \Om} (\rr) \right| < \eps \right\}  \label{eq:Oeps}
\end{align}
is in the tubular neighborhood of $\partial \Om$ of radius $\dbar$,
namely in
$\Odbar \assign \left\{ \rr \in \RR^2 \hspace{0.17em} \stt \left|
    d_{\partial \Om} (\rr) \right| < \dbar \right\}$. For any
$0 < \eps < \dbar$ we also set
$\Od \assign \left\{ \rr \in \RR^2 \stt 0 \leqslant d_{\partial \Om}
  (\rr) < \eps \right\}$. Since $\Om$ is a simply connected $C^2$
domain, there exists a $C^1$ projection map
$\bpi : \Odbar \rightarrow \partial \Om$ such that
$\rr = \bpi (\rr) + d_{\partial \Om} (\rr) \nnn (\bpi (\rr))$ for
every $\rr \in \Odbar$ and
\begin{equation}
  \grad d_{\partial \Om} (\rr) = \nnn (\bpi (\rr)) .
\end{equation}
In particular, $\left| \nabla d_{\partial \Om} (\rr) \right| = 1$ for
every $x \in \Odbar$, and the values of $\bpi$ may be parametrized by
the arclength of $\partial \Omega$. In what follows, we always assume
that $0 < \eps < \dbar$ so that the tubular neighborhood theorem
holds.

    
For any $0 < \eps < \dbar$ we set $\Om_{\eps} \assign \Om \cup \Od$,
which represents the domain in the plane occupied by the ferromagnetic
film, and consider the family of cutoff
functions
\begin{equation}
  \eta_{\eps} (\rr) \assign \eta \left( \frac{d_{\partial \Om}
  (\rr)}{\eps} \right) \label{mdelta}
\end{equation}
defined by a {\tmem{non-increasing}} function $\eta \in C^{0,1}(\RR)$
such that
\begin{equation}
  \eta (t) \equiv 1 \quad \text{for $t \in (- \infty, 0)$,} \quad \eta (t)
  \equiv 0 \quad \text{for $t \in (1, + \infty)$} \label{eq:expreta} .
\end{equation}
We further assume that on $(-\infty,1]$ the function $\eta$ is
continuously differentiable, but allow $\eta'(t)$ to jump at $t = 1$
in accordance with the behavior of $|\overline \m|$ at the film edge
observed in section \ref{sec:micr-deriv-reduc}.


Note that for every $\eps > 0$ we have $\eta_{\eps} (\rr) \equiv 1$
whenever $d_{\partial \Om} (\rr) \leqslant 0$ and
$\eta_{\eps} (\rr) \equiv 0$ when
$d_{\partial \Om} (\rr) \geqslant \eps$. In other words, $\eta_{\eps}$
is a cutoff function whose support is included in the closure of
$\Om_{\eps}$ (i.e.,
$\supp (\eta_{\eps}) \subseteq \overline \Om_{\eps}$) such that
$\eta_{\eps} \equiv 1$ on $\Om$. We observe the following identities:
\begin{equation}
  \nabla \eta_{\eps} (\rr) = \frac{1}{\eps} \eta' \left( \frac{d_{\partial
  \Om} (\rr)}{\eps} \right)  \nnn (\bpi (\rr)) = - \left| \nabla
  \eta_{\eps} (\rr) \right|  \nnn (\bpi (\rr)) . \label{eq:derud}
\end{equation}
In particular, $\grad \eta_{\eps} \in L^\infty ( \R^2; \R^2)$ and
$\supp (\grad \eta_{\eps}) \subseteq \overline \Om_{\eps} \setminus
\Om = \overline{\Od}$.


\subsection{The micromagnetic energy}

Given a configuration
$\m^\eps \in H^1_\mathrm{loc}(\Om_\eps; \Stwo^2)$, we extend it by zero
outside $\Omega_\eps$ to define the magnetization in the whole plane,
so that the physically observable magnetization is
$\m_{\eps} (\rr) = \eta_{\eps} (\rr) \m^{\eps} (\rr)$ for all
$\rr \in \R^2$, after a rescaling of length from \eqref{eq:Omeps}. The
rescaled cutoff length is now $\eps \ll 1$ and the domain
$\Omega_\eps$ has $O(1)$ size, converging to a fixed domain $\Omega$
from outside as $\eps \to 0$. As we will see later, the precise shape
of the cutoff function $\eta$ will prove not to play any role in the
limiting behavior of the energy analyzed in this paper. Further care
is needed, however, to specify a representation of the micromagnetic
energy that is well-suited for the analysis of those limits.

We first need to define a space where the different terms in our
micromagnetic energy are all bounded. To make sure this is the case, we
clearly need that
$\int_{\Om_{\eps}} \eta_{\eps}^2 | \nabla \m^{\eps} |^2 \dd \rr <
\infty$. Therefore, we assume that $\m^{\eps} \in \myspace$ where
$\myspace$ stands for the weighted Sobolev (metric) space defined by
\begin{equation}
  \myspace \assign \left\{ \m^{\eps} \in L^2(\RR^2; \R^3) :
  \m^{\eps}_{| \Om_{\eps}} \in H^1_{\tmop{loc}} ( \Om_{\eps},
  \Stwo^2 ), \eta_{\eps} \grad \m^{\eps}_{| \Om_{\eps}} \in L^2
  ( \Om_{\eps}), \m^{\eps} \equiv 0 \text{ in } \RR^2
  \setminus \Om_{\eps} \right\} .
\end{equation}
Note that elements of $\myspace$ are identically zero outside of
$\Om_{\eps}$. We view $\myspace$ as a metric subspace of $\myspacefull,$
where
\begin{multline}
  \qquad \qquad \myspacefull \assign \big\{ \tmmathbf{u}^{\eps} \in
  L^2 (\RR^2; \R^3) \cap L^{\infty} (\RR^2; \R^3) : \\
  \tmmathbf{u}^{\eps}_{| \Om_{\eps}} \in H^1_{\tmop{loc}} (
  \Om_{\eps}; \R^3), \eta_{\eps} \grad \tmmathbf{u}^{\eps}_{|
    \Om_{\eps}} \in L^2 ( \Om_{\eps}; \R^6), \tmmathbf{u}^{\eps}
  \equiv 0 \text{ in } \RR^2 \setminus \Om_{\eps} \big\} . \qquad
  \qquad
\end{multline}
We can similarly introduce the ``limit'' space
\begin{equation}
  H_0(\R^2; \Stwo^2) \assign \left\{ \m \in L^2(\RR^2; \R^3) :
  \m_{| \Om} \in H^1 ( \Om,
  \Stwo^2 ), \m \equiv 0 \text{ in } \RR^2
  \setminus \Om \right\} .
\end{equation}

We next introduce the notation $\m^{\eps} = (\msq, \mpa^{\eps} )$,
where $\msq \in \R^2$ and $\mpa^{\eps} \in \R$ give the components of
$\m^\eps$ that are perpendicular and parallel to the material easy
axis $\pm \tmmathbf{e}_3$, respectively. The nonlocal contribution
from the stray field energy can then be seen to be proportional to
\begin{equation}
  \mathcal{W}_{\eps} (\m^{\eps}) = \frac{1}{2 | \ln \eps |} \Vdelta
  ( \eta_{\eps} \msq) - \frac{1}{2 | \ln \eps |}
  \Vtildedelta ( \eta_{\eps} \mpa^{\eps}), \label{Wbulkdelta}
\end{equation}
defined for every $\m^{\eps} \in \myspace$, with
\begin{align}
  \Vdelta ( \eta_{\eps} \msq) & \assign   \int_{\R^2}
  \int_{\R^2} \frac{\divv ( \eta_{\eps} \msq) ( \rr )
  \divv ( \eta_{\eps} \msq) ( \sss)}{| \rr - \sss
  |}  \hspace{0.17em} \dd \rr  \dd \sss,  \label{eq:exprmag1}\\
  \Vtildedelta ( \eta_{\eps} \mpa^{\eps}) & \assign 
  \int_{\R^2} \int_{\R^2} \frac{\nabla \left( \eta_{\eps} \mpa^{\eps}
  \right) ( \rr ) \cdot \nabla \left( \eta_{\eps} \mpa^{\eps}
  \right) ( \sss) }{| \rr - \sss | }  \hspace{0.17em} \dd \rr  \dd
  \sss .  \label{eq:exprmag2}
\end{align}
The equivalence of the above expression with the one appearing in
\eqref{eq:Eeps} for smooth functions can be seen via the Fourier
representation.  Note that both $\Vdelta$ and $\Vtildedelta$ are
nonnegative, a result that can be easily shown in the Fourier domain
by Parseval--Plancherel identity.

The DMI contribution to the energy is proportional to
\begin{equation}
  \DMIdelta (\m^{\eps}) \assign \int_{\R^2} \eta_{\eps }^2 \left(
    \mpa^{\eps}  \divv  
    \mpr^{\eps} - \mpr^{\eps} \cdot \nabla \mpa^{\eps} \right)
  \dd \rr .
\end{equation}
The remaining terms may be defined analogously.  Note that the space
$\myspace$ depends on $\eps$ and, therefore, is not well-suited for
$\Gamma$-convergence arguments. Therefore, since
$\m_{\eps} \in L^2 (\R^2; \R^3)$, we can use a penalization to
formulate the $\Gamma$-convergence results in $L^2 (\R^2; \R^3)$ with
the agreement that the energy is infinite outside $\myspace$.
Furthermore, the space $L^2 (\R^2; \R^3)$ also provides a natural
topology for the $\Gamma$-convergence.

The simplified micromagnetic energy defined for
$\m^\eps \in L^2 (\R^2; \R^3)$ that disregards the anisotropy and the
Zeeman terms takes the form
\begin{equation}
  \mathcal{G}_{\eps} ( \m^{\eps}) \assign \left\{
    \begin{array}{lll}  \displaystyle
      \int_{\RR^2} \eta_{\eps}^2 \left| \grad \m^{\eps} \right|^2  \dd \rr +
      \lambda \DMIdelta (\m^{\eps}) +\gamma_\eps
      \mathcal{W}_{\eps} 
      (\m^{\eps}) & 
                    \text{if } \m^{\eps} \in \myspace,\\
      + \infty &   \text{otherwise} ,
  \end{array} \right. \label{eq:enGdelta}
\end{equation}
where the precise dependence of $\gamma_\eps$ on $\eps$ will be
specified in the sequel for various $\Gamma$-limits. As is common in
the studies of $\Gamma$-convergence, we will simply write
$\m^\eps \to \m$ as $\eps \to 0$, always tacitly implying that for any
sequence of $\eps_n \to 0$ we have $\m^{\eps_n} \to \m$ as
$n \to \infty$.

\begin{remark}
  We note that the functional $\mathcal G_\eps$ in \eqref{eq:enGdelta}
  provides a mathematically suitable formulation for the main terms in
  the energy $E_\eps$ defined in \eqref{eq:Eeps} in the considered
  topology of $\Gamma$-convergence. In particular, we have
  $\mathcal G_\eps(\m^{\eps}) = E_\eps(\m^{\eps})$ whenever
  $\m^{\eps} \in \myspace$ with $\alpha = 0$ and $\beta = 0$. Hence
  throughout the rest of the paper we formulate our results in terms
  of the functional $\mathcal G_\eps$ and note that the anisotropy and
  the Zeeman terms are continuous perturbations to $\mathcal G_\eps$
  under $L^2$ convergence and hence can be trivially included in the
  statements of the theorems.
\end{remark}

\subsection{Main results}In this section, we formulate the main results
of the paper. We split our results into four theorems corresponding to
different magnetic regimes previously studied for ferromagnets with
strong in-plane anisotropy.

The first regime we consider is the analogue of the Gioia and James
regime {\cite{Gioia1997}}. We have the following theorem.

\begin{theorem}[Free boundary] \label{t:GJ} Let
  $\mathcal{G}_{\eps} (\m^{\eps})$ be defined on $L^2 (\R^2; \R^3)$ by
  {\tmem{\eqref{eq:enGdelta}}} and let $\gamma_\eps \rightarrow 0$ as
  $\eps \to 0$. We define $\mathcal{G}_0( \m )$ on $L^2 (\R^2; \R^3)$
  by
  \begin{equation}
    \mathcal{G}_0 \left( \m \right) \assign \left\{ \begin{array}{lll} \displaystyle
      \int_{\Om} \left| \grad \m \right|^2 \dd \rr + \lambda \int_{\Om} \left(
      \mpa  \divv  \mpr - \mpr \cdot \nabla \mpa \right)  \dd \rr &  &
      \text{if } \m \in H_0 ( \R^2; \Stwo^2),\\
      + \infty &  & \text{otherwise} .
    \end{array} \right. \label{eq:enG0delta}
  \end{equation}
  Then the following statements hold:
  \begin{enumerateroman}
  \item {\tmem{(Compactness)}} If
    $\limsup_{\eps \to 0} \mathcal{G}_{\eps} (\m^{\eps}) < + \infty$
    then $\m^{\eps} \rightarrow \m$ strongly in $L^2( \RR^2; \RR^3)$
    and $\m^{\eps} \rightharpoonup \m$ weakly in $H^1( \Om; \Stwo^2)$
    for some $\m \in H_0(\R^2; \Stwo^2)$ as $\eps \rightarrow 0$
    (possibly up to a subsequence).
    
  \item {\tmem{(}}{\tmem{$\Gamma \textrm{\text{-}liminf}$}}
    {\tmem{inequality)}} Let $\m^{\eps} \in \myspace$ be such that
    $\m^{\eps} \rightarrow \m$ for some $\m \in H_0(\R^2; \Stwo^2)$
    strongly in $L^2( \RR^2; \RR^3)$ as $\eps \rightarrow 0$. Then
    \begin{align} \liminf_{\eps \rightarrow 0} \mathcal{G}_{\eps} ( \m^{\eps})
      \geqslant \mathcal{G}_0 \left( \m \right) . \end{align}
  \item {\tmem{(}}{\tmem{$\Gamma \textrm{\text{-}limsup}$}}
    {\tmem{inequality)}} Let $\m \in H_0( \R^2; \Stwo^2 )$. Then there
    exists $\m^{\eps} \in \myspace$ such that
    $\m^{\eps} \rightarrow \m$ strongly in $L^2 ( \RR^2; \RR^3)$ as
    $\eps \rightarrow 0$ and
    \begin{align} \limsup_{\eps \rightarrow 0} \mathcal{G}_{\eps} ( \m^{\eps}
       ) \leqslant \mathcal{G}_0 \left( \m  \right) . \end{align}
  \end{enumerateroman}
\end{theorem}

The second regime we study corresponds to the result of Kohn and
Slastikov {\cite{kohn05arma}}, where in the limit of small thickness,
the magnetization prefers to stay in-plane, and a local boundary
contribution corresponding to shape anisotropy replaces the nonlocal
magnetostatic energy. The limiting behavior we obtain for materials
with perpendicular anisotropy is contained in the following
theorem. Here and everywhere below, the values of
$\m \in H_0(\R^2; \Stwo^2)$ on $\partial \Omega$ are understood in the
sense of trace of the Sobolev function
$\m_{|\Omega} \in H^1(\Omega; \Stwo^2)$.

\begin{theorem}[Boundary penalty]
  \label{thm:thmKS}Let $\mathcal{G}_{\eps} ( \m^{\eps})$ be defined on
  $L^2( \RR^2; \R^3)$ by {\tmem{\eqref{eq:enGdelta}}} and
  $\gamma_\eps \to \gamma$ for some $\gamma > 0$ as $\eps \to 0$. We
  define $\mathcal{G}_0^\gamma( \m )$ on $L^2 ( \RR^2; \R^3)$ by
  \begin{equation}
    \mathcal{G}_0^\gamma ( \m ) \assign \left\{ \begin{array}{lll}
      \begin{array}{l}  \displaystyle
        \int_{\Om} \left| \grad \m \right|^2  \dd \rr + \lambda \int_{\Om}
        \left( \mpa  \divv  \mpr - \mpr \cdot \nabla \mpa \right)  \dd \rr\\
        \qquad \qquad \qquad \qquad \qquad + \displaystyle \gamma
        \int_{\partial \Om} \left( 
        (\mpr \cdot \nnn)^2 - \mpa^2 \right)  \dd \mathcal H^1(x)
      \end{array} &  & \text{if } \m \in H_0( \R^2; \Stwo^2),\\
      + \infty &  & \text{otherwise} .
    \end{array} \right. \label{eq:enG0delta1}
  \end{equation}
  Then the following statements hold:
  \begin{enumerateroman}
  \item {\tmem{(Compactness)}} If
    $\limsup_{\eps \to 0} \mathcal{G}_{\eps} ( \m^{\eps}) < +\infty$,
    then $\m^{\eps} \rightarrow \m$ strongly in $L^2( \RR^2; \RR^3)$
    and $\m^{\eps} \rightharpoonup \m$ weakly in $H^1( \Om; \Stwo^2)$
    for some $\m \in H_0(\R^2; \Stwo^2)$ as $\eps \rightarrow 0$
    (possibly up to a subsequence).
    
  \item {\tmem{($\Gamma \textrm{\text{-}liminf}$ inequality)}} Let
    $\m^{\eps} \in \myspace$ be such that $\m^{\eps} \rightarrow \m$
    strongly in $L^2 ( \RR^2; \RR^3 )$ for some
    $\m \in H_0(\R^2; \Stwo^2)$ as $\eps \rightarrow 0$. Then
    \begin{align} \liminf_{\eps \rightarrow 0} \mathcal{G}_{\eps} ( \m^{\eps} )
      \geqslant \mathcal{G}_0^\gamma \left( \m \right) . \end{align}
   \item {\tmem{(}}{\tmem{$\Gamma \textrm{\text{-}limsup}$}}
     {\tmem{inequality)}} Let $\m \in H_0( \R^2; \Stwo^2)$. Then there
     exists $\m^{\eps} \in \myspace$ such that
     $\m^{\eps} \rightarrow \m$ strongly in $L^2( \RR^2; \RR^3)$ as
     $\eps \rightarrow 0$ and
     \begin{align} \limsup_{\eps \rightarrow 0} \mathcal{G}_{\eps} ( \m^{\eps} )
       \leqslant \mathcal{G}_0^\gamma \left( \m \right) . \end{align}
  \end{enumerateroman}
\end{theorem}

The following two results are fundamentally different from what exists
in the micromagnetic literature for in-plane materials. This is due to
the fact that they correspond to magnetic regimes where the shape
anisotropy of the micromagnetic energy is penalized in the limit and,
as a result, the magnetization acquires Dirichlet conditions at the
boundary. For in-plane materials this regime is impossible as it leads
to a singular behavior of the micromagnetic energy due to a
topological obstruction. For materials with perpendicular anisotropy,
however, these regimes provide the micromagnetic energy describing the
behavior of magnetic skyrmions (see \cite{mmss:cmp23}) and therefore
are of utter importance.

We formulate our results in the following theorems corresponding to
local and nonlocal versions of the micromagnetic energies. Note that
everywhere below the statement $\mpa = \pm 1$ on $\partial \Omega$
means that either $\mpa = 1$ on $\partial \Omega$ or $\mpa = -1$ on
$\partial \Omega$, in the sense of trace. In other words, the trace of
$\m$ is constant on $\partial \Omega$ and is equal to either
$\tmmathbf{e}_3$ or $-\tmmathbf{e}_3$.

Our first theorem yields a local limiting micromagnetic energy with
the Dirichlet boundary condition.

\begin{theorem}[Clamped boundary, local energy]
  \label{thm:locengy}
  Let $\mathcal{G}_{\eps} (\m^{\eps})$ be defined on
  $L^2 \left( \RR^2; \R^3 \right)$ by {\tmem{\eqref{eq:enGdelta}}}
  with $\gamma_\eps \rightarrow \infty$ and
  ${\gamma_\eps | \ln \eps |^{-1}} \to 0$ as $\eps \to 0$. We define
  $\tilde{\mathcal{G}}_0 \left( \m \right)$ on
  $L^2\left( \RR^2 ; \R^3 \right)$ by
  \begin{equation}
    \tilde{\mathcal{G}}_0 \left( \m \right) \assign \left\{ \begin{array}{ll} \displaystyle
      \int_{\Om} \left| \grad \m \right|^2  \dd \rr + \lambda \int_{\Om}
      \left( \mpa  \divv  \mpr - \mpr \cdot \nabla \mpa \right)  \dd \rr & 
      \begin{array}{l} 
        \text{if } \m \in H_0 \left( \R^2; \Stwo^2 \right)  \\ \text{and } \mpa =
        \pm 1 \tmop{on} \partial \Om, \\ ~
      \end{array}\\
      + \infty &  ~\, \text{otherwise} .
    \end{array} \right. \label{eq:enG0delta2}
  \end{equation}
  Then the following statements hold:
  \begin{enumerateroman}
  \item {\tmem{(Compactness)}} If
    $\limsup_{\eps \to 0} \left( \mathcal{G}_{\eps} \left( \m^{\eps}
      \right) + \gamma_\eps \mathcal{H}^1 \left( \partial \Om \right)
    \right) < +\infty$ then $\m^{\eps} \rightarrow \m$ strongly in
    $L^2 ( \RR^2; \RR^3)$ and $\m^{\eps} \rightharpoonup \m$ weakly in
    $H^1 \left( \Om; \Stwo^2 \right)$ for some
    $\m \in H_0(\R^2; \Stwo^2)$ with $\mpa =
        \pm 1$ on
    $\partial \Om$ as $\eps \rightarrow 0$ (possibly up to a
    subsequence).
    
  \item {\tmem{($\Gamma \textrm{\text{-}liminf}$ inequality)}} Let
    $\m^{\eps} \in \myspace$ be such that $\m^{\eps} \rightarrow \m $
    strongly in $L^2 \left( \RR^2; \RR^3 \right)$ for some
    $\m \in H_0(\R^2; \Stwo^2)$ with $\mpa =
        \pm 1 $ on
    $\partial \Om$ as $\eps \rightarrow 0$. Then
    \begin{align} \liminf_{\eps \rightarrow 0} \left( \mathcal{G}_{\eps} \left(
          \m^{\eps} \right) + \gamma_\eps \mathcal{H}^1 (\partial
        \Omega) \right) \geqslant \tilde{\mathcal{G}}_0 \left( \m
      \right) . \end{align}
  \item {\tmem{($\Gamma \text{-} \mathrm{limsup}$ inequality)}} Let
    $\m \in H_0 \left( \R^2; \Stwo^2 \right)$ with
    $\mpa =
        \pm 1$ on $\partial \Om$. Then there exists
    $\m^{\eps} \in \myspace$ such that $\m^{\eps} \rightarrow \m$
    strongly in $L^2 \left( \RR^2; \RR^3 \right)$ as
    $\eps \rightarrow 0$ and
    \begin{align} \limsup_{\eps \rightarrow 0} \left( \mathcal{G}_{\eps} \left(
          \m^{\eps} \right) + \gamma_\eps \mathcal{H}^1 (\partial
        \Omega) \right) \leqslant \tilde{\mathcal{G}}_0 \left( \m
      \right) . \end{align}
  \end{enumerateroman}
\end{theorem}

The next theorem provides a new type of a reduced nonlocal
micromagnetic energy for yet stronger dipolar interaction. To state
the theorem, we need to introduce some additional notation. This is
due to the fact that in the considered regime the stray field of the
film edge continues to contribute to the limit energy and, therefore,
needs to be properly accounted for. We define the quantity
\begin{equation}
  \Ddelta \assign \int_{\R^2} \int_{\R^2} \frac{\nabla
    \eta_{\eps} (\rr) \cdot \nabla \eta_{\eps} (\sss) }{| \rr - \sss |}
  \dd \rr  \dd \sss, \label{eq:Ddeltadef}
\end{equation}
which will be shown to give, up to the factor of $\frac12 \nu$, the
leading order behavior of the energy $\mathcal G_\eps$. In fact, this
constant can be seen to be the energy of the ferromagnetic state
$\m = \pm \tmmathbf{e}_3 \chi_\Om$ and to the leading order satisfies
$D_\eps = 2 |\ln \eps| \, \mathcal H^1(\partial \Omega) + O(1)$ when
$\eps \to 0$, as can be seen from Lemma \ref{lem:Deps}.

\begin{theorem}[Clamped boundary, nonlocal energy]
  \label{thm:nonlocengy}Let $\mathcal{G}_{\eps} (\m^{\eps})$ be
  defined on $L^2\left( \RR^2; \R^3 \right)$ by
  {\tmem{\eqref{eq:enGdelta}}} with $\gamma_\eps = \nu | \ln \eps |$
  for some $\nu > 0$. We define
  $\tilde{\mathcal{G}}_0^\nu \left( \m \right)$ on
  $L^2 \left( \RR^2 ; \R^3 \right)$ by
  \begin{equation}
    \tilde{\mathcal{G}}_0^\nu  \left( \m \right) \assign \left\{ \begin{array}{ll}
      \begin{array}{l}  \displaystyle
        \int_{\Om} \left| \grad \m \right|^2  \dd \rr  + \nu \int_\Om
        \tmmathbf{b} \cdot \nabla \mpa \, \dd \rr \\
        \quad \displaystyle + \lambda \int_{\Om}
        \left( \mpa  \divv  \mpr - \mpr \cdot \nabla \mpa \right)  \dd
        \rr \\ 
        \displaystyle \quad + \frac{\nu}{2}  \int_{\Om} \int_{\Om} \frac{\divv  \mpr(
        \rr )  \divv  \mpr (\sss)}{| \rr - \sss |}  \hspace{0.17em} \dd
        \rr  \dd \sss\\
        \displaystyle  \quad - \frac{\nu}{2}  \int_\Om
        \int_\Om \frac{\nabla  
        \mpa ( \rr ) \cdot \nabla  \mpa
         ( \sss)}{| \rr - \sss | }  \dd \rr  \dd
        \sss
      \end{array} &  \begin{array}{l}
        \text{if } \m \in H_0 \left( \R^2; \Stwo^2 \right)\\
        \text{and } \mpa = \pm 1 \tmop{on} \partial \Om,
      \end{array}\\
      + \infty &  \text{otherwise} ,
    \end{array} \right. \label{eq:enG0delta3}
\end{equation}
where $\tmmathbf{b}$ is defined in \eqref{eq:b}. Then the following
statements hold:
  \begin{enumerateroman}
  \item {\tmem{(Compactness)}} If
    $\limsup_{\eps \to 0} \left( \mathcal{G}_{\eps} \left( \m^{\eps}
      \right) + \frac{\nu}{2} \Ddelta \right) < +\infty$ then
    $\m^{\eps} \rightarrow \m$ strongly in
    $L^2 \left( \RR^2; \RR^3 \right)$ and
    $\m^{\eps} \rightharpoonup \m$ weakly in
    $H^1 \left( \Om; \Stwo^2 \right)$ for some
    $\m \in H_0(\R^2; \Stwo^2)$ with $\mpa = \pm 1$ on
    $\partial \Omega$ as $\eps \rightarrow 0$ (possibly up to a
    subsequence), where $D_\eps$ is defined in \eqref{eq:Ddeltadef}.
    
  \item {\tmem{($\Gamma \textrm{\text{-}liminf}$ inequality)}} Let
    $\m^{\eps} \in \myspace$ satisfy $\m^{\eps} \rightarrow \m$
    strongly in $L^2 \left( \RR^2; \RR^3 \right)$ for some
    $m \in H_0(\R^2; \Stwo^2)$ with
    $\mpa =
        \pm 1$ on $\partial \Omega$ as
    $\eps \rightarrow 0$, then
    \begin{align} \liminf_{\eps \rightarrow 0}  \left( \mathcal{G}_{\eps} \left(
       \m^{\eps}  \right) + \frac{\nu}{2} \Ddelta \right)
       \geqslant \tilde{\mathcal{G}}_0^\nu  \left( \m  \right) . \end{align}
   \item {\tmem{($\Gamma \text{-} \mathrm{limsup}$ inequality)}} Let
     $\m \in H_0 \left( \Om; \Stwo^2 \right)$ with
     $\mpa =
        \pm 1$ on $\partial \Omega$. Then
     there exists $\m^{\eps} \in \myspace$ that satisfies
     $\m^{\eps} \rightarrow \m$ strongly in
     $L^2 \left( \RR^2; \RR^3 \right)$ as $\eps \rightarrow 0$ and
     \begin{align} \limsup_{\eps \rightarrow 0} \left( \mathcal{G}_{\eps} \left(
           \m^{\eps} \right) + \frac{\nu}{2} \Ddelta \right) \leqslant
       \tilde{\mathcal{G}}_0^\nu \left( \m \right) . \end{align}
  \end{enumerateroman}
\end{theorem}

The above results provide us with four comprehensive two-dimensional
micromagnetic models to study the magnetization behavior in ultrathin
films with perpendicular magnetic anisotropy and Dzyaloshinskii-Moriya
interaction.  It is also clear that the above results can be
supplemented with anisotropy and Zeeman energies as those are just
continuous perturbations. In particular, this provides a rigorous
justification to the formal limit statements in section
\ref{sec:model}.

\section{Auxiliary lemmas}\label{sec:techlemmas}

Throughout the rest of the paper, unless stated otherwise all the
constants in the statements and the proofs depend only on $\Omega$ and
$\eta$.
We begin by providing several important technical lemmas. The first
two lemmas concern pointwise estimates for the singular integral
\begin{equation}
  f_{\eps} (\rr) \assign
  \int_{\Od} \frac{\left| \grad \eta_{\eps} (\sss) \right|}{| \sss - \rr |} 
  \hspace{0.17em} \dd \sss. \label{eq:fdeltanew0}
\end{equation}
We start with an estimate which for every $0 < \eps < \dbar$ gives a
precise control on $f_{\eps} (x)$ when $x \in \Odpm$, where $\Odpm$ is
defined in \eqref{eq:Oeps}.

\myprop{\begin{lemma}
    \label{Lemma:LemmaforK1} There exists $\dbar, C > 0$ such that for
    any $0 < \eps < \dbar$ and $\rr \in \Odpm$ there holds
  \begin{equation}
    2 | \ln \eps | - C \leqslant f_{\eps} ( \rr ) \leqslant 2
    | \ln \eps | + C. \label{eq:boundfdxoverlog}
  \end{equation}
\end{lemma}}

\begin{proof}
  Recalling the notations of section \ref{sec:geom}, for any
  $|t| < \dbar$, the curve
  \begin{equation}
    \varphi_t : s \in I_{\partial \Om} \mapsto \varphi (s) + t \nnn (s) ,
    \label{eq:paramgammatnew}
  \end{equation}
  is a parameterization of the set
  $\partial \Om_t \assign \left\{ \sigma \in \Odpepsbar \of d_{\partial
      \Om} (\sigma) = t \right\}$. Since $\Om$ is of class $C^2$,
  there exists $\bar{\delta} > 0$ sufficiently small such that for any
  $0 < \delta < \bar{\delta}$ there exists a curvature-dependent
  constant $0 < \alpha (\delta) < \delta$ for which there holds
  \begin{equation}
    \mathcal{H}^1 ( \bpi ( \Odpepsbar \cap B_{\alpha (\delta)}
    ( \rr ) )) < \delta \quad \forall x \in
    \Odpepsbar . \label{eq:alphaeps}
  \end{equation}
  In particular, there holds that
  $\mathcal{H}^1 ( \bpi ( \partial \Om_t \cap B_{\alpha
        (\delta)} ( \rr ) ) ) < \delta$ for any
  $|t|< \dbar$, and we can always assume that $\dbar$ and
  $\bar{\delta}$ are tuned sufficiently small, so that for any
  $|t|< \dbar$ and any $0 < \delta < \bar{\delta}$, the set
  $\bpi ( \partial \Om_t \cap B_{\alpha (\delta)} ( \rr )
  )$ is connected. In this way, the arc
  $\bpi ( \partial \Om_t \cap B_{\alpha (\delta)} ( \rr )
  )$ can be parameterized through the restriction of $\varphi$
  to a suitable subinterval of $I_{\partial \Om}$. In other words, we
  assume that $\delta$ is sufficiently small so that for any
  $|t| < \dbar$ and any $\rr = \varphi_t (s_0) \in \Odpepsbar$ one has
  \begin{equation}
    \partial \Om_t \cap B_{\alpha (\delta)} ( \rr ) \subseteq
    \varphi_t \left( I_{\delta} ( \rr ) \right), \quad
    I_{\delta} ( \rr ) \assign [s_0 - \delta / 2, s_0 +
    \delta / 2] . \label{eq:parcovercurvenew}
  \end{equation}
  Now, let $\eps < \dbar$. For what follows, it is convenient to set
  \begin{equation}
    \kappa_{\partial \Om} \assign \sup \left\{ | \kappa (\sigma) | \stt \sigma
    \in \partial \Om \right\}, \label{eq:kdeltabar}
  \end{equation}
  where $\kappa (\sigma)$ stands for the curvature of $\partial \Om$
  at the point $\sigma \in \partial \Om$
  ({\tmabbr{cf.}}~\eqref{eq:Frenet1}). Clearly, for any
  $\rr \in \Odpm$ and any
  $\sss \in \Od \setminus B_{\alpha (\delta)} \left( \rr \right)$ one
  has $\left| \sss - \rr \right| \geqslant \alpha (\delta)$.  We
  denote by $S_{\alpha (\varepsilon)}(x)$ the small sector
  around $x \in \Odpm$ defined by
  \begin{align}
    S_{\alpha (\delta)} ( \rr ) \assign \left\{ \varphi_t \left(
    I_{\delta} ( \rr ) \right) \right\}_{|t| <
    \eps}.     
  \end{align}
  Clearly,
  $S_{\alpha (\delta)} ( \rr ) \supseteq B_{\alpha (\delta)} ( \rr )
  \cap \Od$, therefore we can decompose $f_{\eps} (\rr)$ as
  \begin{align}
    f_{\eps} (\rr)
    & 
        = \psi_{\eps}(x) +  g_{\partial \Om}^\eps ( \rr ) 
        ,  \label{eq:fdeltanew} 
  \end{align}
  where
  \begin{equation}
    \psi_{\eps} (x) \assign \frac{1}{\eps} \int_{\Od \cap S_{\alpha
    (\delta)} ( \rr )} \frac{\left| \eta' ( d_{\partial
    \Om} (\sss) / \eps ) \right|}{| \sss - \rr |}  \hspace{0.17em} \dd
    \sss
  \end{equation}
  and by the coarea formula the remainder term
  $g_{\partial \Om}^\eps ( \rr )$ satisfies the uniform bound
  \begin{equation}
    \left| g_{\partial \Om}^\eps ( \rr ) \right| \leqslant
    \frac{1}{\alpha (\delta)} \int_0^1 \int_{\partial \Om} | \eta' (t)
    | \cdot
    | 1 + \eps t \kappa (\sigma) | 
    \hspace{0.17em} \dd \sigma \dd t \; \leqslant c_{\delta}, \quad
    c_{\delta} \assign \frac{\left( 1 + \bar\varepsilon
        \kappa_{\partial \Om} \right) \mathcal{H}^1 \left( \partial \Om
      \right)}{\alpha (\delta)} . \label{eq:gdOm}
  \end{equation}

  It remains to estimate $\psi_{\eps} (x)$. For that we observe that
  with $x = \sigma + \eps t\tmmathbf{n} (\sigma)$, $|t| < 1$, and
  by the coarea formula we have
  \begin{eqnarray}
    \psi_{\eps} (x)
    & \eqs &
               \int_0^1 \int_{\varphi \left( I_{\delta} ( \rr )
               \right)} | \eta' (r) | \frac{1 + \eps r \kappa
               (\omega)  }{| \sigma +
               \eps t\tmmathbf{n} (\sigma) - \omega - \eps r\tmmathbf{n} (\omega) |} 
               \dd \mathcal H^1(\omega) \dd r. 
  \end{eqnarray}
  We estimate the denominator of the integrand as
  \begin{align}
    \omega - \sigma + \eps (r\tmmathbf{n} (\omega) - t\tmmathbf{n} (\sigma))
    & \eqs  \varphi (s) - \varphi  (s_0) + \eps (r\tmmathbf{n} (s) -
             t\tmmathbf{n} (s_0)) \nonumber\\
    & \eqs  (s - s_0) \tmmathbf{\tau} (s_0) + \eps (r - t) \tmmathbf{n}
    (s_0)  \nonumber \\
   & \qquad  \qquad \qquad \qquad  + \eps r (\tmmathbf{n} (s) -\tmmathbf{n} (s_0)) + O (|s -
             s_0|^2), 
    \label{eq:2bcomb1}
  \end{align}
  where
  $| O (|s - s_0|^2) | \leqslant \kappa_{\partial \Om} (s - s_0)^2
  $. On the other hand, since $0 < r < 1$ and $|t| < 1$, and
  $\tmmathbf{n}(s)$ is Lipschitz continuous with constant
  $\kappa_{\partial \Omega}$, we also have
  \begin{align}
    | \eps r (\tmmathbf{n} (s) -\tmmathbf{n} (s_0)) + O (|s - s_0|^2) |
    & \leqslant 
                 \kappa_{\partial \Omega} (\eps + | s - s_0 |) | s - s_0 |
                 \nonumber\\
    & \leqslant  \kappa_{\partial \Omega} (\eps + | s - s_0 |) | (s - s_0)
                  \tmmathbf{\tau} (s_0) + \eps (r - t) \tmmathbf{n}
                  (s_0) | . 
    \label{eq:2bcomb2}
  \end{align}
  Hence, combining \eqref{eq:2bcomb1} and \eqref{eq:2bcomb2} we obtain
  that for $\bar \eps, \delta < \frac12 \kappa_{\partial \Om}^{-1}$
  and $|s - s_0| < \delta$ we have
  \begin{align}
    | \omega - \sigma + \eps (r\tmmathbf{n} (\omega) - t\tmmathbf{n}
    (\sigma)) | & \leqslant  (1 + \kappa_{\partial \Omega} \eps
                              + \kappa_{\partial 
                              \Omega} | s - s_0 |) | (s - s_0) \tmmathbf{\tau} (s_0) + \eps (r - t)
                              \tmmathbf{n} (s_0) |,  \\
    | \omega - \sigma + \eps (r\tmmathbf{n} (\omega) - t\tmmathbf{n}
    (\sigma)) | & \geqslant  (1 - \kappa_{\partial \Omega} \eps
                              - \kappa_{\partial 
                              \Omega} | s - s_0 |) | (s - s_0) \tmmathbf{\tau} (s_0) + \eps (r - t)
                              \tmmathbf{n} (s_0) | . 
  \end{align}
  Overall, we get that for any
  $x = \sigma + \eps t\tmmathbf{n} (\sigma)$, $|t| < 1$, there
  holds
  \begin{equation}
    \frac{1 - \eps \kappa_{\partial \Om}}{1 + \eps \kappa_{\partial \Om}}
    \int_0^1 | \eta' (r) | \varrho_{\eps}^+ \left(
      \frac{\eps^2}{\delta^2} (r - t)^2 \right)  \dd r \leqslant
    \psi_{\eps} (x) \leqslant \frac{1 + \eps \kappa_{\partial \Om}}{1 -
      \eps \kappa_{\partial \Om}} \int_0^1 | \eta' (r) | \varrho_{\eps}^- \left(
      \frac{\eps^2}{\delta^2} (r - t)^2 \right)  \dd r,
    \label{eq:estlogd1}
  \end{equation}
  with
  \begin{align}
    \varrho_{\eps}^{\pm} \left( \frac{\eps^2}{\delta^2} (r - t)^2
    \right) & \assign  \int_{I_{\delta} ( \rr )} \frac{1}{(1
    \pm a^{\pm}_{\eps} | s - s_0 |) | (s - s_0) \tmmathbf{\tau} (s_0) +
    \eps (r - t) \tmmathbf{n} (s_0) |} \dd s \\
    & \eqs  \int_{s_0 - \delta / 2}^{s_0 + \delta / 2} 
    \frac{1}{(1 \pm a^{\pm}_{\eps} | s - s_0 |) \sqrt{(s - s_0)^2 + \eps^2
    (r - t)^2}}  \dd s \\
    & \eqs  2 \int_0^{\frac{1}{2}} \frac{1}{\left( 1 \pm a^{\pm}_{\eps}
    \delta s \right) \sqrt{s^2 + \frac{\eps^2}{\delta^2} (r - t)^2}} \dd
    s,  \label{eq:exprrho1}
  \end{align}
  where $a_{\eps}^{\pm} = \frac{\kappa_{\partial \Omega}}{1 \pm \eps
  \kappa_{\partial \Omega}}>0$.
  
A direct integration yields
  \begin{equation}
    \varrho_{\eps}^{\pm} (\beta) = \rho(\beta) \mp
    \int_0^{\frac12} \frac{2 a^\pm_\eps \delta s}{ (1 \pm a^\pm_\eps
      \delta s) \sqrt{s^2 + \beta}} \dd s,
    \label{eq:exprrho2}
  \end{equation}
  where
    \begin{equation}
    \varrho (\beta) \assign \ln \left( \frac{1 + 2
        \beta + \sqrt{1 + 4 
          \beta}}{2 \beta} \right).
    \label{eq:exprrho2lemma8}
  \end{equation}
  By inspection, for all $\beta < 1$ we have
  \begin{equation}
    - \ln \beta - C \leqslant \varrho_{\eps}^\pm
    (\beta) \leqslant - \ln \beta + C, 
    \label{eq:estlogd2}
  \end{equation}
  for some $C > 0$ universal.  Also note that 
  \begin{eqnarray}
    \int_0^1 | \eta' (r) | \cdot | \ln (r - t)^2 |  \dd r  \leq C \|
    \eta' \|_\infty, 
    \label{eq:estlogd3}
  \end{eqnarray}
  for some $C > 0$ universal.
  Therefore, given \eqref{eq:estlogd2} and \eqref{eq:estlogd3}, from
  the relation \eqref{eq:estlogd1} we infer the estimate
  \begin{equation}
    2 | \ln \eps | - C \leqslant \psi_{\eps} ( \rr )
    \leqslant 2 | \ln \eps | + C, \label{eq:estlogd1new}
  \end{equation}
  for some $C > 0$, assuming that $\eps < \delta$ and that
  $\delta < \bar \delta$ is chosen sufficiently small.  The conclusion
  of the lemma then follows from the previous estimate,
  \eqref{eq:fdeltanew}, and \eqref{eq:gdOm}.
\end{proof}

\begin{remark}
  We point out that Lemma \ref{Lemma:LemmaforK1} and all the
  subsequent results, which rely on Lemma \ref{Lemma:LemmaforK1},
  remain valid if $\eta(t)$ is only H\"older continuous at $t = 1$
  (with an arbitrary H\"older exponent).
\end{remark}

We next prove a pointwise bound on the function $f_\eps (x)$ outside
$\Odpm$.

{\myprop{\begin{lemma}
      \label{1stlemma0} There exist $\dbar, C > 0$ depending only on
      $\Om$, such that for any $0 < \eps < \dbar$ and
      $x \in \Omega \backslash \Odpm$ there holds
\begin{align}
  0 \leqslant f_\eps(x)
  \leqslant C \big( 1 + |\ln (\mathrm{dist}(x,
  \partial \Omega))| \big).
    \end{align}
\end{lemma}}}

\begin{proof}
  As in the proof of Lemma \ref{Lemma:LemmaforK1}, using the coarea
  formula we infer that for any $\rr \in \Omega \backslash \Odpm$
  and sufficiently small $\eps$ there holds
  \begin{equation}
    f_{\eps} ( \rr )
    \eqs \int_{\partial \Om} \int_0^1 
             | \eta' (t) | 
             \frac{ 1 + \eps t \kappa (\sigma)  }{\left| \rr - \sigma - \eps t
             \nnn (\sigma) \right|} 
             \tmop{dt} \dd \mathcal H^1( \sigma ), 
  \end{equation}
  where $\kappa(\sigma)$ is a curvature at point
  $\sigma \in \partial\Omega$.  Furthermore, we claim that
  $|x - \sigma - \eps t \nnn(\sigma)| \geq \frac12 |x - \sigma|$ for
  every $x \in \Omega \backslash \Odpm$, $\sigma \in \partial \Om$ and
  $t \in [0,1]$, provided that $\eps$ is sufficiently small. To see
  this, notice that the estimate trivially holds when
  $|x - \sigma| \geq 2 \eps$ or when
  $(x - \sigma) \cdot \nnn(\sigma) < 0$. At the same time, in the
  opposite case we can estimate
  $|x - \sigma| \leq 2 |(x - \sigma) \cdot \tmmathbf{\tau}(\sigma)|
  \leq 2 |x - \sigma - \eps t \nnn(\sigma)|$ for all $\eps$
  sufficiently small in view of the regularity of $\partial
  \Omega$. Thus, we have
    \begin{align}
       \label{eq:fepsdomlndist}
      0 \leqslant f_\eps(x) 
      \leqslant 4 \int_{\partial \Omega} \frac{1}{ |x - \sigma|} \dd
      \mathcal H^1(\sigma)  \leqslant C \big( 1 + |\ln (\mathrm{dist}(x,
      \partial \Omega))| \big),
    \end{align}
    for some $C > 0$ and all $\eps$ small enough.
    \end{proof}

    As an immediate consequence of Lemma~\ref{Lemma:LemmaforK1} and
    Lemma~\ref{1stlemma0} we have the following result.

{\myprop{\begin{lemma}
      \label{1stlemma} There exist $\dbar, C > 0$ such that for every
      $0 < \eps < \dbar$ and every $x \in \Om_\eps$ there holds
  \begin{equation}
      |f_\eps(x)|
      \leqslant C \big( 1 + |\ln (\eps + \mathrm{dist}(x,
      \partial \Omega))| \big).
  \end{equation}
  In particular, for every $p\geq 1$ there is $C_p > 0$ such that
  $\|f_{\eps} \|_{L^p \left( \Om_\eps \right)} \leqslant C_p$.
\end{lemma}}}

We will also need a sharp estimate for the quantity $D_\eps$
introduced in \eqref{eq:Ddeltadef}.

{\myprop{\begin{lemma}
      \label{lem:Deps} There exist $\dbar, C > 0$ such that for every
      $0 < \eps < \dbar$ there holds
  \begin{eqnarray}
  2 |\ln \eps| \, \mathcal H^1(\partial \Omega) - C \leqslant D_\eps
    \leqslant 2 |\ln \eps| \, \mathcal H^1(\partial \Omega) + C.
  \end{eqnarray}
\end{lemma}}}

\begin{proof}
  Observe that for all $\dbar$ sufficiently small we have
  \begin{align}
    D_\eps= \int_{\Od} \int_{\Od} \frac{|\nabla \eta_\eps(x)| \, |\nabla
    \eta_\eps(y)|}{|x - y|} \, \dd x \dd y - \frac12
    \int_{\Od} \int_{\Od} |\nabla \eta_\eps(x)| \, |\nabla
    \eta_\eps(y)|  \frac{|\nnn (x) - \nnn(y)|^2}{ |x - y|} \, \dd x
    \dd y, \label{Deps2ints}
  \end{align}
  with an abuse of notation $\nnn(x) \assign \nnn(\bpi(x))$ for
  $x \in \Od$. Notice that by the $C^2$ regularity of
  $\partial \Omega$ the second integral in \eqref{Deps2ints} is
  uniformly bounded when $\eps \to 0$, as $\nnn(x)$ is Lipschitz
  continuous. At the same time, by the coarea formula and Lemma
  \ref{Lemma:LemmaforK1} the first integral in \eqref{Deps2ints} is
  \begin{align}
    \int_{\Od} \int_{\Od}  & \frac{|\nabla \eta_\eps(x)| \, |\nabla
    \eta_\eps(y)|}{|x - y|} \, \dd x \dd y 
    = \int_{\Od}
    |\nabla \eta_\eps(x)| f_\eps(x) \, \dd x \notag \\
    & = \int_0^1
    \int_{\partial \Om} |\eta'(t)| \, f_\eps(\sigma + \eps t \nnn(\sigma)
    ) (1 + \eps t \kappa(\sigma))  \, \dd \mathcal H^1(\sigma) \,
    \dd t   
     = 2 |\ln \eps| \, \mathcal H^1(\partial \Om) + O(1),
  \end{align}
  as $\eps \to 0$, which implies the statement of the lemma.
\end{proof}

The next key technical lemma provides a comparison of an integral
involving $\nabla \eta_\eps$ tested against a bounded Sobolev function
with that of the same integral evaluated on the trace of that Sobolev
function on $\partial \Omega$.

{\myprop{\begin{lemma}
      \label{prop:boundsforj2j3}
      There exist constants $\dbar, c > 0$ such that for every
      $0 < \mu < 1$ there holds
  \begin{equation}
    \frac{1}{c} \int_{\Od} \left| \nabla \eta_{\eps} (\rr)
    \right|  |u (\rr) - u (\bpi (\rr)) |  \hspace{0.17em} \dd \rr  \leqslant 
    \mu \| u \|_{L^{\infty} (\Od)} + | \ln \mu |  \|
    \eta_{\eps} \nabla u \|_{L^1 (\Od)}, 
    \label{eq:estmperpuniflemma}
  \end{equation}
  for any
  $u \in L^{\infty} (\Od) \cap W^{1, 1}_{\tmop{loc}} ( \Om_{\eps})$.
\end{lemma}}}

\begin{proof}
  Again, with the usual abuse of notation
  $\tmmathbf{n} (x) \assign \tmmathbf{n} (\bpi (x))$, we have
  $\rr = \bpi (\rr) + d_{\partial \Om} (\rr) \nnn ( \rr )$ for every
  $x \in \Od$. Moreover, as noted in \eqref{eq:derud}, we have that
  \begin{equation}
    \left| \nabla \eta_{\eps} (\rr) \right| =  \frac{1}{\eps}  \left| \eta' (
      d_{\partial \Om} (\rr) / \eps ) \right|.
    \label{eq:tempnprimedelta0}
  \end{equation}
  Therefore, taking the precise representative of $u$, by its
  differentiability on the lines $\bpi(x) = y$ for $\mathcal
  H^1$-a.e. $y \in \partial \Omega$ and monotonicity of $\eta(t)$ we
  have
  \begin{align}
    \left| \left( u (\rr) - u (\bpi (\rr)) \right) \nabla \eta_{\eps} (\rr)
    \right| & \eqs  \left| u (\bpi (\rr) + d_{\partial \Om} (\rr) \nnn (
    \rr )) - u (\bpi (\rr)) \right|  \left| \nabla \eta_{\eps} (\rr)
    \right| \nonumber\\
    & \leqslant  \left| \nabla \eta_{\eps} (\rr) \right| \cdot
    \int_0^{d_{\partial \Om} (\rr)} \left| \partial_t \left[ u \left( \bpi
    (\rr) + t \nnn ( \rr ) \right) \right] \right|  \dd t
    \nonumber\\
            & \leqslant  \frac{\left| \eta' \left( d_{\partial
                          \Om} (\rr) / \eps
                          \right) \right|}{\eps \eta \left( d_{\partial
                          \Om} (\rr) / \eps \right) }
                          \int_0^{d_{\partial \Om} (\rr)} \eta (t
                          / \eps) \left| \grad u \left( 
                          \bpi (\rr) + t \nnn ( \rr ) \right) \right|  \dd t. 
    \label{eq:tempnprimedelta}
  \end{align}
  
  For $0 < \lambda < 1$ we decompose $\Od$ as
  $\Od =\mathcal{O}_{\lambda \eps}^+ \cup \left( \Od
    \backslash\mathcal{O}_{\lambda \eps}^+ \right)$ and focus
  separately on $\mathcal{O}_{\lambda \eps}^+$ and
  $\Od \backslash\mathcal{O}_{\lambda \eps}^+$.  Starting with
  $\mathcal{O}_{\lambda \eps}^+$, we observe that both $\bpi(x)$ and
  $\nnn(x)$ are constant along the normal direction, so by
  \eqref{eq:tempnprimedelta0}, \eqref{eq:tempnprimedelta} and the
  coarea formula we infer that
  \begin{align}
    \int_{\mathcal{O}_{\lambda \eps}^+}
    & \left| \nabla
      \eta_{\eps} (\rr)
      \right| \left| u (\rr) - u (\bpi (\rr)) \right|  \dd x \leqslant 
      \int_{\mathcal{O}_{\lambda \eps}^+} \frac{\left| \eta' \left( d_{\partial
      \Om} (\rr) / \eps
      \right) \right|}{\eps \eta \left( d_{\partial
      \Om} (\rr) / \eps \right) } \left( \int_0^{\lambda \eps} \eta (t
      / \eps)
      \left| \, \grad u \left( \bpi (\rr) + t \nnn ( \rr )_{_{_{_{}}}}
      \right)_{_{}} \right| \dd t \right) \dd x \nonumber\\
    & \eqs   \int_0^{\lambda} \int_{\partial \Om} \frac{|
      \eta' (h) |}{\eta (h)} \left( \int_0^{\lambda \eps}
      \eta (t
      / \eps) \left| \, \grad u \left( \sigma + t \nnn
      (\sigma)_{_{_{_{}}}} \right)_{_{}} \right| \dd t \right)  | 1 +
      \eps h \kappa
      (\sigma) |  \dd \mathcal H^1 (\sigma) \dd h \nonumber\\
    & \leqslant   3  |\ln \eta(\lambda)|
      \int_0^{\lambda \eps} \int_{\partial \Om}
      \eta (t
      / \eps) \left| \, \grad u \left( \bpi (\rr) + t \nnn \left( \rr
      \right)_{_{_{_{}}}} \right)_{_{}} \right| | 1 + t \kappa (\sigma) | \dd
      \mathcal H^1 (\sigma) \dd t,  \label{eq:tempnprimedelta2}
  \end{align}
  provided $\dbar < \frac12 \kappa_{\partial \Om}^{-1} $.
  Thus
  \begin{multline}
    \int_{\mathcal{O}_{\lambda \eps}^+} \left| \nabla \eta_{\eps}
      (\rr) \right| \left| u (\rr) - u (\bpi (\rr)) \right| \dd x
    \\
    \leqslant 3 | \ln \eta (\lambda) | \int_0^{\lambda \eps}
    \int_{\partial \Om} \eta_{\eps} (t) \left| \, \grad u \left( \bpi
        (\rr) + t \nnn ( \rr )_{_{_{_{}}}} \right)_{_{}} \right| | 1 +
    t \kappa (\sigma) | \dd \mathcal H^1(\sigma) \dd t
    \\
    \eqs 3 | \ln \eta (\lambda) | \cdot \| \eta_{\eps} \nabla u
    \|_{L^1 (\mathcal{O}_{\lambda \eps}^+)} .
    \label{eq:firstestimuminuproj}
  \end{multline}

  On the other hand, for the part on
  $\Od \backslash\mathcal{O}_{\lambda \eps}^+$, we have, again by
  \eqref{eq:tempnprimedelta0} and coarea formula, that
  \begin{eqnarray}
    \int_{\Od \backslash\mathcal{O}_{\lambda \eps}^+} \left| \nabla
    \eta_{\eps} (\rr) \right|  |u (\rr) - u (\bpi (\rr)) |  \hspace{0.17em}
    \dd \rr & \eqs & \frac{1}{\eps} \int_{\Od
                     \backslash\mathcal{O}_{\lambda 
                     \eps}^+}  \left| \eta' \left( d_{\partial \Om} (\rr) / \eps \right)
                     \right|  |u (\rr) - u (\bpi (\rr)) |
                     \hspace{0.17em} \dd \rr \nonumber\\ 
            & \leqslant & 4 |\partial \Om|  \left(
                          \int_{\lambda}^1 | \eta' (t) 
                          |  \dd t \right) \| u \|_{L^{\infty} (\Od)} \nonumber\\
            & = & 4 |\partial \Om|  \eta (\lambda)  \| u
                  \|_{L^{\infty} (\Od)} .  \label{eq:secondestimuminuproj}
  \end{eqnarray}
  Overall, combining the estimates \eqref{eq:firstestimuminuproj} and
  \eqref{eq:secondestimuminuproj}, we get that
  \begin{eqnarray}
    \frac{1}{c} \int_{\Od} \left| \nabla \eta_{\eps} (\rr)
    \right|  |u (\rr) - u (\bpi (\rr)) |  \hspace{0.17em} \dd \rr
    & \leqslant &
                  \eta (\lambda) \| u \|_{L^{\infty} (\Od)} + | \ln \eta
                  (\lambda) | \, \| \eta_{\eps} \nabla u \|_{L^1 (\Od)} ,
  \end{eqnarray}
  for some $c > 0$.  The previous estimate holds for every
  $0 < \lambda < 1$. Since $\eta (\lambda)$ maps $[0, 1]$ surjectively
  onto $[0, 1]$, setting $\lambda \assign \eta^{- 1} (\mu)$ we get
  that for every $\mu \in (0, 1)$ there holds:
  \begin{eqnarray}
    \frac{1}{c} \int_{\Od} \left| \nabla \eta_{\eps} (\rr)
    \right|  |u (\rr) - u (\bpi (\rr)) |  \hspace{0.17em} \dd \rr
    & \leqslant &
                  \mu \| u \|_{L^{\infty} (\Od)} + | \ln \mu |  \|
                  \eta_{\eps} \nabla u \|_{L^1 (\Od)},
    \label{eq:estimatoimprove} 
  \end{eqnarray}
  which proves the $L^1$-estimate \eqref{eq:estmperpuniflemma}.
\end{proof}

Our next lemma gives a bound that will be useful to estimate the
interior contribution of the bulk charges to the micromagnetic
energy. Note that for $\uuu \in H^1(\Om; \R^2)$ a straightforward
interpolation estimate between the $\mathring{H}^{-1/2}$ norm of
$\divv \uuu$ and the $L^2$ norms of $\uuu$ and $\nabla \uuu$ would
have held true if $\uuu$ vanished at the boundary of
$\Omega$. However, the presence of a non-zero trace on
$\partial \Omega$ requires some additional care due to a logarithmic
failure of this interpolation. A counterexample to the latter is
provided by $\uuu$, which is equal to the outward unit normal at the
projection point to the boundary of $\Omega$ multiplied by a cutoff
function making $\uuu$ zero at distances greater than $\delta$ from
$\partial \Omega$ (for a related phenomenon, see \cite{desimone06}).
  
\myprop{\begin{lemma}
    \label{Lemma:divdivu} There exists a constant $C > 0$ depending
    only on $\Om$ such that for any $\delta \in (0, \frac12)$ and any
    $\uuu \in H^1(\Omega; \R^2)$ there holds
  \begin{align}
    \int_{\Om} \int_{\Om} \frac{\divv \uuu ( \rr ) \divv \uuu ( \sss
      )}{\left| \rr - \sss \right|} \dd \rr \dd \sss
    \leqslant \delta \| \nabla \uuu \|_{L^2(\Omega)}^2 + C \delta^{-1} \| \uuu
    \|_{L^2(\Omega)}^2 + C |\ln \delta| \, \| \uuu \|_{L^2(\partial
      \Omega)}^2 . \label{eq:divL2s}
  \end{align}
\end{lemma}}

\begin{proof}
  For a nonnegative cutoff function $\omega \in C^\infty(\R)$ satisfying
  $\omega(t) = 1$ for all $t \leq 1$ and $\omega(t) = 0$ for all
  $t \geq 2$, we write
  \begin{align}
    \frac{1}{|x - y|} = G_\delta(x - y) + H_\delta(x - y),
  \end{align}
  where $H_\delta(x) := |x|^{-1} \omega(c |x| / \delta)$ with $c > 0$
  to be fixed shortly and
  $G_\delta(x) := |x|^{-1} (1 - \omega(c |x| / \delta))$. For the
  contribution of $H_\delta$, by Young's inequality for convolutions
  we have
  \begin{align}
    \label{eq:Hdel}
    \left| \int_{\Om} \int_{\Om}
    H_\delta(x - y) \divv  \uuu ( \rr )  \divv  \uuu (
    \sss ) \dd \rr \dd \sss \right| \leqslant 2 \| H_\delta
    \|_{L^1(\R^2)} \| \nabla \uuu \|_{L^2(\Om)}^2
    \leqslant \delta \| \nabla \uuu \|_{L^2(\Om)}^2,
  \end{align}
  for a suitable choice of $c > 0$ depending only on $\omega$.  For
  the contribution of the even function $G_\delta$, we integrate by
  parts and apply the Cauchy-Schwarz inequality to obtain
  \begin{align}
    \label{eq:Gdel}
    \int_{\Om} \int_{\Om}
    &
      G_\delta(x - y) \divv \uuu ( \rr ) \divv
      \uuu (
      \sss ) \dd \rr \dd \sss  = \int_{\Om} \int_{\Om} \uuu ( \rr )
      \cdot \nabla^2_{xy} G_\delta(x - 
      y) \uuu ( \sss ) \dd \rr \dd \sss \notag \\
    & - 2 \int_{\Om}
      \int_{\partial \Om} (\uuu ( y ) \cdot \nnn(y)) \, \uuu(x) \cdot
      \nabla_x G_\delta(x - y) \, \dd \mathcal H^1(y) \dd x \notag \\
    & + \int_{\partial \Om}
      \int_{\partial \Om} (\uuu ( x ) \cdot \nnn(x)) \,  (\uuu ( y )
      \cdot \nnn(y)) \, G_\delta(x - y) \, \dd \mathcal H^1(x) \, \dd
      \mathcal H^1(y) 
      \notag \\
    & \leqslant \max_{y \in \Om} \| \nabla^2 G_\delta(\cdot - y)
      \|_{L^1(\Om)} \| \uuu \|_{L^2(\Om)}^2 +
      \max_{y \in \partial \Om} \| G_\delta(\cdot - y)
      \|_{L^1(\partial \Omega)} \| \uuu \|_{L^2(\partial \Om)}^2  
      \notag \\ 
    & + 2 \max_{y \in \partial \Omega} \| \nabla G_\delta(\cdot
      - y) \|_{L^1(\Omega)}^{1/2} \max_{y \in \Omega} \| \nabla G_\delta(\cdot
      - y) \|_{L^1(\partial \Omega)}^{1/2} \| \uuu
      \|_{L^2(\Om)} \| \uuu \|_{L^2(\partial \Om)} \notag \\
    & \leqslant C \delta^{-1} \| \uuu \|_{L^2(\Om)}^2 + C |\ln
      \delta| \, \| \uuu \|_{L^2(\partial \Om)}^2 + C \delta^{-1/2}
      |\ln \delta|^{1/2} \| \uuu \|_{L^2(\Om)} \, \| \uuu
      \|_{L^2(\partial \Om)}, 
  \end{align}
  for some $C > 0$ depending only on $\Om$.  The conclusion follows by
  combining \eqref{eq:Hdel} and \eqref{eq:Gdel}, after an application
  of Young's inequality.
\end{proof}

As a corollary to this lemma, we have the following result for vector
fields that are uniformly bounded. 

\begin{lemma}
  \label{lm:con}
  Let $\delta > 0$, and for $\eps > 0$ let
  $\uuu_{\eps} \in H^1 ( \Om; \R^2 )$ with
  $\left| \uuu_{\eps} \right| \leqslant 1$ in $\Om$. Then there exists
  a constant $C > 0$ depending only on $\Om$ such that
  \begin{equation}
    \int_{\Om} \int_{\Om}
    \frac{\divv  \uuu_{\eps} ( \rr )  \divv  \uuu_{\eps} (
      \sss )}{\left| \rr - \sss \right|} \dd \rr \dd \sss \leqslant \delta
    \left\| \nabla \uuu_{\eps} \right\|^2_{L^2 \left( \Om  \right)} +
    \frac{C}{\delta} . \label{eq:inequ}
  \end{equation}
  Moreover, if $\uuu_{\eps} \rightharpoonup \uuu$ weakly in
  $H^1 \left( \Om \right)$ as $\eps \to 0$ then
  \begin{equation}
    \int_{\Om} \int_{\Om}
    \frac{\divv  \uuu_{\eps} ( \rr )  \divv  \uuu_{\eps} (
      \sss )}{\left| \rr - \sss \right|} \dd \rr \dd \sss \xrightarrow{\eps
      \rightarrow 0} \int_{\Om} \int_{\Om}
    \frac{\divv  \uuu ( \rr )  \divv  \uuu (
      \sss )}{\left| \rr - \sss \right|} \dd \rr \dd \sss \label{eq:divconv}
  \end{equation}
  and
  \begin{equation}
    \int_{\partial \Om} \int_{\Om}
    \frac{\divv  \uuu_{\eps} ( \rr )  (  \uuu_{\eps} (
      \sss ) \cdot \nnn(y))}{\left| \rr - \sss \right|} \, \dd \rr \, \dd
    \mathcal 
    H^1(\sss) \xrightarrow{\eps 
        \rightarrow 0} \int_{\partial \Om} \int_{\Om}
    \frac{\divv  \uuu ( \rr )  (  \uuu (
      \sss ) \cdot \nnn(y))}{\left| \rr - \sss \right|} \, \dd \rr \, \dd
    \mathcal  
    H^1(\sss) \label{eq:divconvtrace}    
  \end{equation}
\end{lemma}

\begin{proof}
  The estimate in \eqref{eq:inequ} is an immediate corollary to Lemma
  \ref{Lemma:divdivu}. To prove \eqref{eq:divconv}, we first note that
  since $|x|^{-1} * \divv \uuu$ belongs to $L^2(\Om)$ by Young's
  inequality for convolutions, it is enough to prove that the
  left-hand side of \eqref{eq:divconv} goes to zero when
  $\uuu_\eps \rightharpoonup 0$ weakly in $H^1(\Om)$ as $\eps \to
  0$. However, the latter is true by \eqref{eq:divL2s} in view of
  boundedness of $\| \nabla \uuu \|_{L^2(\Om)}$, strong convergence of
  $\uuu$ to zero in $L^2(\Om)$ and $L^2(\partial \Om)$ by compact
  Sobolev and trace embeddings, and arbitrariness of $\delta >
  0$. Similarly, since the integral over $\partial \Om$ in the
  right-hand side of \eqref{eq:divconvtrace} defines a function of $x$
  that belongs to $L^2(\Omega)$ by the last inequality in
  \eqref{eq:fepsdomlndist}, it is enough to show
  \eqref{eq:divconvtrace} when $\uuu_\eps \rightharpoonup 0$ weakly in
  $H^1(\Om)$ as $\eps \to 0$. The latter follows via an application of
  the Cauchy-Schwarz inequality:
    \begin{align}
      \left| \int_{\partial \Om} \int_{\Om}
      \frac{\divv  \uuu_{\eps} ( \rr )  (  \uuu_{\eps} (
      \sss ) \cdot \nnn(y))}{\left| \rr - \sss \right|} \, \dd \rr \, \dd
      \mathcal 
      H^1(\sss) \right| \leqslant \left( \int_{\partial \Om} \int_{\Om}
      \frac{ (  \uuu_{\eps} (
      \sss ) \cdot \nnn(y))^2}{\left| \rr - \sss \right|^{3/2}} \, \dd
      \rr \, \dd \mathcal H^1(\sss) \right)^{1/2} \notag \\
      \times \left( \int_{\partial \Om} \int_{\Om}
      \frac{ |\divv  \uuu_{\eps} ( \rr ) |^2}{\left| \rr - \sss \right|^{1/2}}
      \, \dd \rr \, \dd     \mathcal 
      H^1(\sss) \right)^{1/2} \leqslant C \| \nabla \uuu_\eps 
      \|_{L^2(\Om)} \| \uuu_\eps \cdot \nnn \|_{L^2(\partial \Om)} \to 0, 
    \end{align}
    as $\eps \to 0$, by the compact trace embedding, where $C > 0$
    depends only on $\Om$.  
\end{proof}

\section{Analysis of the magnetostatic energy} \label{sec:asympt}

In this section we carry out an analysis of the magnetostatic energy
which contains two propositions describing the behavior of the
nonlocal terms in the stray field energy $\mathcal{W}_{\eps}$
({\tmabbr{cf.}}~\eqref{Wbulkdelta}), corresponding to the in-plane and
out-of-plane magnetization components. We will then use these results
to prove our main theorems formulated in section \ref{sec:discussion}.

We start by proving the following proposition for the nonlocal term
due to the in-plane component of the magnetization. Note that here we
need a stronger result than that of the type proved for the
three-dimensional micromagnetic energy by Kohn and Slastikov in
\cite{kohn05arma} in order to go beyond the regime studied there (see
Theorems \ref{thm:locengy} and \ref{thm:nonlocengy}).

\begin{proposition}
  \label{th:mainresultdemag} There exist $\dbar, C > 0$ such that if
  $0 < \eps < \dbar$ and $\m^\eps \in \myspace$, then the
  magnetostatic energy for the in-plane component
  {\tmem{(}}{\tmabbr{cf.}}~{\tmem{\eqref{eq:exprmag1})}}
  \begin{equation}
    \Vdelta ( \eta_{\eps} \msq ) =
    \int_{\R^2} \int_{\R^2} \frac{\divv ( \eta_{\eps} \msq)
      ( \rr ) \, \divv ( \eta_{\eps} \msq ) ( \sss)}{| \rr - \sss |}
    \hspace{0.17em} \dd \rr  \dd \sss
  \end{equation}
  satisfies
  \begin{equation}
    \begin{array}{l}
      \left| \mathcal{V} \left( \eta_{\eps}  \msq \right) -\mathcal{V}_{\Om
      \times \Om} (\msq) + \left. 2\mathcal{V}_{\partial \Om
      \times \Om} (\msq) - 2 | \ln \eps | \, \| \msq \cdot
      \nnn \|^2_{L^2( \partial \Om)} \right| \qquad \qquad
      \right.\\
      \qquad \qquad  \leqslant \;
      C_{\eps} \left( 1 + \| \eta_{\eps} \nabla
      \msq \|^2_{L^2( \Om_{\eps})} \right) + C \| \msq \cdot \nnn
      \|_{L^2(\partial \Om )} \left( 1 + \| \eta_{\eps} \nabla
      \msq \|_{L^2( \Om_{\eps})} \right),
    \end{array} \label{eq:2beproved}
  \end{equation}
  where $C_{\eps} \rightarrow 0$ as $\eps \rightarrow 0$, and
  \begin{align}
    \mathcal{V}_{\Om \times \Om} ( \msq )
    & \assign  \int_{
                \Om} \int_{\Om} \frac{ \divv \msq ( \rr )  \, \divv \msq
                ( \sss )}{| \rr - \sss |}  \dd \sss \,
                \dd \rr, \\
    \mathcal{V}_{\partial \Om \times \Om}( \msq )
    & \assign 
                \int_{\partial \Om} \int_{\Om}  \frac{( \nnn
                \cdot \msq ) (\sigma) \,  \divv \msq ( \sss ) }{|
                \sigma - \sss |}  \dd \sss \, \dd \mathcal H^1(\sigma). 
  \end{align}
\end{proposition}

\begin{proof}
  As is common in the analysis of the limiting behaviors of nonlocal
  energy functionals, we decompose $\mathcal{V}$ into the sum of
  several terms and estimate each term separately.
  
  First, expanding the divergence and exploiting the symmetry in the
  $\rr, \sss$ variables, we write
  $\mathcal{V} \left( \uuu_{\eps} \right) \backassign I_1 + 2 I_2 +
  I_3$, where
  \begin{align}
    I_1 & \assign  \int_{\R^2} \int_{\R^2} \frac{\grad \eta_{\eps} (\rr)
    \cdot \msq ( \rr )  \hspace{0.17em} \grad \eta_{\eps} (\sss)
    \cdot \msq (\sss)}{\left| \rr - \sss \right|}  \dd \rr \dd \sss, \\
    I_2 & \assign   \int_{\R^2} \int_{\R^2} \frac{\nabla \eta_{\eps} (\sss)
    \cdot \msq (\sss) \eta_{\eps} (\rr)  \divv  \msq (\rr)}{\left| \rr -
    \sss \right|}  \hspace{0.17em} \dd \rr \dd \sss, \\
    I_3 & \assign  \int_{\R^2} \int_{\R^2} \frac{\eta_{\eps} (\rr)  \divv
    \msq ( \rr ) \eta_{\eps} ( \sss)  \divv \msq
    ( \sss)}{\left| \rr - \sss \right|}  \hspace{0.17em} \dd \rr
    \dd \sss . 
  \end{align}
  Note that $I_1, I_2, I_3$ depend on $\eps$ but we suppress this for
  ease of notation. We then have
  \begin{align}
    \text{left-hand side of \eqref{eq:2beproved}} \leqslant & \left|\, I_3
    -\mathcal{V}_{\Om \times \Om} (\msq) \, \right| + 2 \left|\,
    \mathcal{V}_{\partial \Om \times \Om} (\msq) + I_2 \,\right| \\
    & \qquad \qquad \qquad \qquad
    + \left| \, I_1 - 2 | \ln \eps | \| \msq \cdot \nnn \, \|^2_{L^2
    \left( \partial \Om \right)} \right|,  \label{eq:firstsubd}
  \end{align}
  and we proceed by estimating the terms on the right-hand side of the
  previous relation.
  
  \vspace{8pt}
  {\noindent}{\tmname{Step 1.}}  {\noindent}{\tmname{Estimate of
      $I_3 -\mathcal{V}_{\Om \times \Om} \left( \msq \right)$.}} We
  split this term as
  $I_3 -\mathcal{V}_{\Om \times \Om} \left( \msq \right) \eqs 2 L_2 +
  L_3$, where
  \begin{equation}
    L_2 \assign \int_{\Om} \int_{\Od} \frac{\eta_{\eps} (\rr)
      \divv \msq 
      ( \rr )  \divv \msq ( \sss)}{\left| \rr - \sss
      \right|}  \hspace{0.17em} \dd \rr \dd \sss,
  \end{equation}
  \begin{equation}
    L_3 \assign \int_{\Od} \int_{\Od} \frac{\eta_{\eps} (\rr)  \divv \msq
      ( \rr ) \eta_{\eps} ( \sss)  \divv \msq \left(
        \sss \right)}{\left| \rr - \sss \right|}  \hspace{0.17em} \dd
    \rr \dd \sss  .
  \end{equation}
  To estimate $L_2$ and $L_3$, we use the Young's inequality for
  convolutions
  \begin{equation}
    \int_{\RR^2} \left| f ( \rr ) \right|  \left| (g \ast K) \left(
        \rr \right) \right|  \dd \rr  \; \leqslant \; \| f \|_p \| g \|_s \| K
    \|_r\label{eq:GYI}
  \end{equation}
  with $p = r = \frac43$ and $s = 2$. Observing that since $\Om$ is
  bounded, there exists a ball $U$ centered at the origin such that
  $\rr - \sss \in U$ for every $\rr, \sss \in \overline{\Om}_{\eps}$
  and $\| | \cdot |^{- 1} \|_{L^r (U)} \leqslant C$ we obtain that
  $L_2 \leqslant C \| \eta_{\eps} \divv \msq \|_{L^p( \Od)} \|
  \eta_{\eps} \divv \msq \|_{L^2 ( \Om_{\eps}),}$. But by H\"older's
  inequality there holds
  $\| \eta_{\eps} \nabla \msq \|_{L^p (\Od)} \leqslant \| \eta_{\eps}
  \nabla \msq \|_{L^2( \Om_{\eps})} |\Od|^{1/4} \leqslant C \eps^{1 /
    4} \| \eta_{\eps} \nabla \msq \|_{L^2( \Om_{\eps})},$ for some
  constant $C > 0$, provided $\eps$ is small enough. Hence
  \begin{equation}
    L_2 \leqslant C \eps^{1 / 4} \| \eta_{\eps} \nabla \msq
    \|^2_{L^2(\Om_{\eps})} \label{esL2}.
  \end{equation}
  The very same estimate is true for $| L_3 |$ and, therefore, we
  conclude that for every $\eps$ sufficiently small there holds
  \begin{equation}
    \left| I_3 -\mathcal{V}_{\Om \times \Om} (\msq) \right|
    \leqslant 2 | L_2 | + | L_3 | \leqslant C \eps^{1 / 4}  \left\|
    \eta_{\eps} \nabla \msq \right\|^2_{L^2 ( \Om_{\eps})},
    \label{estI3n}
  \end{equation}
  for some positive constant $C > 0$ and all $\eps$ small enough.
  
    \vspace{8pt}
  {\noindent}{\tmname{Step 2. Estimate of
      $\mathcal{V}_{\partial \Om \times \Om} (\msq) + I_2$.}} Our aim
  here is to show that
  \begin{equation}
    \left| \mathcal{V}_{\partial \Om \times \Om} (\msq) + I_2
    \right| \leqslant C_\eps \left( 1 + \left\| \eta_{\eps} \nabla
        \msq \right\|^2_{L^2 ( \Om_{\eps})} \right),
    \label{eq:VOmcrossOmplusI2} 
  \end{equation}
  for some $C_\eps > 0$ such that $C_\eps \to 0$ as $\eps \to 0$.  To
  estimate $\mathcal{V}_{\partial \Om \times \Om} (\msq) + I_2$, we
  use the fact that
  $\tmop{supp}_{\RR^2} \grad \eta_{\eps} \subseteq \overline{\Od}$ and
  $\left| \msq \right| \leqslant 1$ in $\Om_{\eps}$ to obtain
  ({\tmabbr{cf.}}~\eqref{eq:derud})
  \begin{align}
    - I_2
    & \eqs  \int_{\Od} \int_{\Om_{\eps}}
             \frac{\left| \nabla \eta_{\eps} (\sss) \right| \left( \nnn \left( \sss
             \right) \cdot \msq (\sss) \right) \eta_{\eps} (\rr) \divv \msq
             (\rr)}{\left| \rr - \sss \right|}  \hspace{0.17em} \dd \rr  \dd \sss \\
    & \backassign  M_1 + M_2 ,
  \end{align}
  with
  \begin{align}
    M_1 & \assign  \int_{\Od} \int_{\Om_{\eps}}
                    \frac{\left| \nabla \eta_{\eps} ( \sss) \right| \left( \nnn
                    \cdot \msq \right) (\bpi ( \sss)) \eta_{\eps} (\rr)  \divv 
                    \msq (\rr)}{\left| \rr - \sss \right|} \dd \rr \dd \sss \\
    M_2 & \assign   \int_{\Od} \int_{\Om_{\eps}}
                    \frac{\left| \nabla \eta_{\eps} ( \sss) \right| \left[
                    \left( \nnn \cdot \msq \right) (\sss) - \left( \nnn \cdot \msq \right)
                    (\bpi ( \sss)) \right] \eta_{\eps} (\rr)  \divv  \msq
                    (\rr)}{\left| \rr - \sss \right|}  \hspace{0.17em} \dd \rr  \dd \sss . 
  \end{align}
  Clearly, we have
  \begin{align} \left| \mathcal{V}_{\partial \Om \times \Om} (\msq) + I_2
     \right| \leqslant \left| \mathcal{V}_{\partial \Om \times \Om} \left(
     \msq \right) - M_1 \right| + | M_2 |, \end{align}
  and we want to show that
  \begin{equation}
    | M_2 | \leqslant C_{\eps} \left( 1 + \left\| \eta_{\eps}  \grad  \msq
    \right\|^2_{L^2 ( \Om_{\eps})} \right) \quad \text{and}
    \quad \left| \mathcal{V}_{\partial \Om \times \Om} (\msq) -
    M_1 \right| \leqslant C_{\eps} \left( 1 + \left\| \eta_{\eps} \nabla
    \msq \right\|^2_{L^2 ( \Om_{\eps})} \right), \label{eqs:M1M2}
  \end{equation}
  with $C_{\eps} \rightarrow 0$ as $\eps \to 0$. To estimate $M_2$, we
  use the Young's inequality for convolutions \eqref{eq:GYI} to obtain
  \begin{equation}
    M_2 \leqslant \| | \cdot |^{- 1} \|_{L^s (U)} \, \left\| | \nabla
      \eta_{\eps} | \left[ \left( \nnn \cdot \msq \right) (\cdot) - \left(
          \nnn \cdot \msq \right) (\bpi (\cdot)) \right] \right\|_{L^p \left( \Od
      \right)} \left\| \eta_{\eps}  \divv  \msq \right\|_{L^2 \left(
        \Om_{\eps} \right)},
  \end{equation}
  for some $p,s \geq 1$ such that
  $\frac{1}{p} + \frac{1}{s} = \frac{3}{2}$. We take
  $p \assign 1 + \alpha$ and
  $s \assign 2 \frac{1 + \alpha}{1 + 3 \alpha}$ with $\alpha > 0$
  sufficiently small so that $1 < s < 2$. Then
  \begin{equation}
    M_2  \leqslant  C \left\| | \nabla \eta_{\eps} | \left[( \nnn
    \cdot \msq ) (\cdot) - ( \nnn \cdot \msq ) (\bpi (\cdot))
    \right]\right\|_{L^{1 + \alpha} (\Od)} \|\eta_{\eps}
    \grad  \msq\|_{L^2 ( \Om_{\eps}),}
  \end{equation}
  with some $C > 0$ depending only on $\Om$ and $s$ such that
  $\| | \cdot |^{- 1} \|_{L^s (U)} \leqslant C$.
  
  We conclude by showing that
  \begin{equation}
    A_{\eps} \assign \left\| | \nabla \eta_{\eps} | \left[ \left( \nnn
    \cdot \msq \right) (\cdot) - \left( \nnn \cdot \msq \right) (\bpi (\cdot))
    \right] \right\|_{L^{1 + \alpha} (\Od)}
    \label{eq:weneedlemmaforthis}
  \end{equation}
  is small, which guarantees that the limit relation in
  \eqref{eqs:M1M2} holds.  Indeed, using interpolation inequality
  $\|f\|_{L^{1 + \alpha}} \, \leqslant \, \|f\|^{\theta}_{L^1}
  \|f\|^{1 - \theta}_{L^{1 + 2 \alpha}}$ with $\alpha > 0$,
  $\theta = \frac{1}{2 (1 + \alpha)}$ and
  $1 - \theta = \frac{1 + 2 \alpha}{2 (1 + \alpha)}$, we immediately
  obtain
  \begin{align} A_{\eps} \; \leqslant \; \left\| | \nabla \eta_{\eps}
      | \left[ \left( \nnn \cdot \msq \right) (\cdot) - \left( \nnn
          \cdot \msq \right) (\bpi (\cdot)) \right]
    \right\|^{\theta}_{L^1 (\Od)} \| 2 \nabla \eta_{\eps} \|^{1 -
      \theta}_{L^{1 + 2 \alpha} (\Od)} . \end{align} Now, recalling
  that $\eta' \in L^{1 + 2 \alpha} (0, 1)$ for some $\alpha$ small
  enough depending on $q$, we get that
  \begin{align}
  \| \nabla \eta_{\eps} \|^{1 - \theta}_{L^{1 + 2 \alpha} (\Od)}
  {\leqslant C \eps^{2 \theta - 1}}.
  \end{align}
   Also, using
  Lemma~\ref{prop:boundsforj2j3} with $\mu = \eps$ and
    Cauchy-Schwarz inequality we obtain
  \begin{align}
    \left\| | \nabla \eta_{\eps} | \left[ \left( \nnn \cdot \msq \right)
    (\cdot) - \left( \nnn \cdot \msq \right) (\bpi (\cdot)) \right]
    \right\|_{L^1 (\Od)} \leqslant C \left( \eps + |\ln \eps| \, |
    \Od |^{1/2}  \| \eta_{\eps} \grad \msq
    \|_{L^2 ( \Om_{\eps})} \right) \notag \\
    \leqslant C' \eps^{1 /
    2} | \ln \eps | \left( 1 + \| \eta_{\eps} \grad \msq
    \|_{L^2 ( \Om_{\eps})} \right), 
  \end{align}
  for some $C, C' > 0$ and all $\eps$ small enough. Now, combining the
  two previous estimates we obtain that
\begin{align} 
  A_{\eps} \leqslant C \eps^{\frac{5 \theta}{2} - 1} | \ln \eps
  |^{\theta} \left( 1 +\| \eta_{\eps} \grad \msq \|_{L^2 (
  \Om_{\eps})} \right)^{\theta} \leqslant C' \eps^{1/8}
  \left( 1 +\| \eta_{\eps} \grad \msq \|_{L^2 ( \Om_{\eps})}
  \right), 
\end{align} 
for $\alpha$ sufficiently small, recalling that $\theta \to \frac12$
as $\alpha \to 0$. Thus by Young's inequality we have
  \begin{align} 
    M_2 \leqslant C \eps^{1 / 8} \left( 1 + \| \eta_{\eps} \grad 
    \msq \|^2_{L^2 ( \Om_{\eps})} \right) ,
\end{align}
which implies the first estimate in \eqref{eqs:M1M2}.
  
It remains to estimate the quantity
$\left| \mathcal{V}_{\partial \Om \times \Om} (\msq) - M_1
\right|$. For that we split it further as $M_1 \backassign N_1 + N_2$,
where
  \begin{align}
    N_1 & \assign  \int_{\Od} \int_{\Om_{\eps}}
          \frac{\left| \nabla 
          \eta_{\eps} ( \sss) \right| \left( \nnn \cdot \msq \right)
          (\bpi ( \sss))  \divv  \msq (\rr)}{\left| \rr - \sss \right|} 
          \dd \rr \dd \sss, \\
    N_2 & \assign  \int_{\Od} \int_{\Od} \frac{\left| \nabla
          \eta_{\eps} ( \sss) \right| \left( \nnn \cdot \msq \right)
          (\bpi ( \sss)) \eta_{\eps} (\rr)  \divv  \msq (\rr)}{\left|
          \rr - \sss \right|}  \dd \rr \dd \sss ,
  \end{align}
  and given that
  $\left| \mathcal{V}_{\partial \Om \times \Om} \left( \msq \right) -
    M_1 \right| \leqslant \left| \mathcal{V}_{\partial \Om \times \Om}
    (\msq) - N_1 \right| + | N_2 |$ we aim at showing that
  \begin{equation}
    \left| \mathcal{V}_{\partial \Om \times \Om} (\msq) - N_1
    \right| \leqslant C_{\eps} \left( 1 + \| \eta_{\eps} \nabla \msq
      \|^2_{L^2 ( \Om_{\eps})} \right) \quad \text{and}
    \quad | N_2 | \leqslant C_\eps \left( 1 +
      \| \eta_{\eps} \nabla \msq \|^2_{L^2( \Om_{\eps}
        )} \right), \label{eqs:N1N2}
  \end{equation}
  with some $C_{\eps} \rightarrow 0$ as $\eps \to 0$. The bound
  \eqref{eqs:N1N2} will prove the second estimate in \eqref{eqs:M1M2}.

  We can estimate $N_2$ as
  \begin{align}
    | N_2 | \leqslant \int_{\Od} f_{\eps} ( \rr ) \,
    \eta_{\eps} (\rr) \left| \divv  \msq (\rr) \right|  \dd
    \rr,
  \end{align}
  where $f_\eps$ is defined in \eqref{eq:fdeltanew0}.  Now we apply
  Lemma~\ref{Lemma:LemmaforK1}, the Cauchy-Schwarz and Young's
  inequalities to obtain that for $\eps$ sufficiently small we have
  \begin{align}
    | N_2 | & \leqslant (2 | \ln \eps | + C) \int_{\Od} 
              \eta_{\eps} (\rr) \, | \divv  \msq (\rr) | \dd
              \rr \nonumber\\
            & \leqslant C'  | \ln \eps | \, | \Od|^{1/2} \|
              \eta_{\eps} \nabla \msq \|_{L^2( \Om_{\eps})}  \notag
    \\ 
            & \leqslant  C''  \eps^{1 / 2} | \ln \eps | \left( 1 + \|
              \eta_{\eps} \nabla \msq \|^2_{L^2( \Om_{\eps})}
              \right) ,  \label{eqN2}
  \end{align}
  for some $C, C', C'' > 0$, yielding the second relation in
  \eqref{eqs:N1N2}.

  To estimate
  $| \mathcal{V}_{\partial \Om \times \Om} ( \msq ) - N_1 |$, we
  observe that
  \begin{equation}
    \left| \mathcal{V}_{\partial \Om \times \Om} ( \msq ) - N_1
    \right|  \eqs \left| \int_{\Om} \divv \msq( \rr )
      \rho_{\eps} ( \rr )  \dd \rr \right| 
  \end{equation}
  with
  \begin{equation}
    \rho_{\eps} ( \rr ) :\eqs \int_{\Od}
    \frac{\left| \nabla \eta_{\eps} ( \sss) \right| \left( \nnn
        \cdot \msq \right) (\bpi ( \sss))}{\left| \rr - \sss \right|} 
    \dd \sss - \int_{\partial \Om}  \frac{ \left( \nnn \cdot \msq
      \right) (\sigma)}{\left| \rr - \sigma \right|}  \dd \mathcal
    H^1(\sigma).
  \end{equation}
  Using the coarea formula, we infer that for any $\rr \in \Om$ there
  holds
  \begin{eqnarray}
    \rho_{\eps} ( \rr )
    & \eqs & \int_{\partial \Om} \int_0^1 
             \left( \nnn \cdot \msq \right) (\sigma) | \eta' (t) | \left(
             \frac{ 1 + \eps t \kappa (\sigma) }{\left| \rr - \sigma - \eps t
             \nnn (\sigma) \right|} - \frac{1}{\left| \rr - \sigma \right|} \right)
             \tmop{dt} \dd \mathcal H^1( \sigma ), 
  \end{eqnarray}
  and by Lebesgue's dominated convergence theorem we have
  $\rho_\eps(x) \to 0$ as $\eps \to 0$ for every $x \in \Om$,
  Furthermore, by $|\msq| \leqslant 1$ and Lemma \ref{1stlemma} we
  have
  \begin{align}
    \label{eq:rhoepsdom}
    |\rho_\eps(x)| \leqslant f_\eps(x) + \int_{\partial \Omega} \frac{1}{|x - \sigma|} \dd 
    \mathcal H^1(\sigma)   
    \leqslant C \big( 1 + |\ln (\mathrm{dist}(x,
    \partial \Omega))| \big),
  \end{align}
  for some $C > 0$ and all $\eps$ small enough.

  From \eqref{eq:rhoepsdom} and the pointwise convergence of
  $\rho_\eps$ to zero as $\eps \to 0$, one can conclude by the
  Lebesgue's dominated convergence theorem that
  $\| \rho_{\eps} \|_{L^2 ( \Om)} \rightarrow 0$ when
  $\eps \rightarrow 0$. Therefore, by the Cauchy-Schwarz and Young's
  inequalities we obtain
  \begin{equation}
    \left| \mathcal{V}_{\partial \Om \times \Om} (\msq) - N_1
    \right| \leqslant C_{\eps} \left( 1 + \| \eta_{\eps} \nabla \msq
    \|^2_{L^2( \Om_{\eps}),} \right)
    \label{eq:missingestimV}
  \end{equation}
  for some $C_{\eps} > 0$ such that $C_{\eps} \rightarrow 0$ when
  $\eps \rightarrow 0$. This concludes the proof of \eqref{eqs:N1N2}
  and, therefore, of \eqref{eq:VOmcrossOmplusI2}.
  
    \vspace{8pt}
  {\noindent}{\tmname{Step 3.}} {\tmname{Estimate of}}
  $\left| I_1 - 2 | \ln \eps | \left\| \msq \cdot \nnn
    \right\|^2_{L^2 \left( \partial \Om \right)} \right|$. By adding
  and subtracting $\msq (\bpi (\rr))$ to $\msq (\rr)$, we rewrite
  $I_1$ as $I_1 = J_1 + 2 J_2 -  J_{3}$, where
  \begin{align}
    J_1 & \assign  \int_{\Od} \int_{\Od} \frac{\nabla \eta_{\eps} (\rr)
    \cdot \msq (\bpi (\rr))  \hspace{0.17em} \nabla \eta_{\eps} (\sss) \cdot
    \msq (\bpi (\sss))}{| \rr - \sss |}  \hspace{0.17em} \dd \rr  \dd \sss, \\
    J_2 & \assign   \int_{\Od} \int_{\Od} \frac{\nabla \eta_{\eps} (\rr)
    \cdot [\msq (\rr) - \msq (\bpi (\rr))] \nabla \eta_{\eps} (\sss) \cdot
    \msq (y)}{| \rr - \sss |}  \hspace{0.17em} \dd \rr  \dd \sss, \\
    J_3 & \assign   \int_{\Od} \int_{\Od} \frac{\nabla \eta_{\eps} (\rr)
    \cdot [\msq (\rr) - \msq (\bpi (\rr))] \nabla \eta_{\eps} (\sss) \cdot
    [\msq (\sss) - \msq (\bpi (\sss))]}{| \rr - \sss |}  \hspace{0.17em} \dd
    \rr  \dd \sss . 
  \end{align}
  In writing the previous relations, we exploited that
  $\supp \nabla \eta_{\eps} \subseteq \overline{\Od}$. Also, to avoid
  cumbersome notations we use the same symbol to denote both $\msq$
  and its trace $\msq_{|\partial \Om}$ on $\partial \Om$. When we
  write $\msq (\bpi (\rr))$ we mean
  $\msq_{| \partial \Om} (\bpi (\rr))$.
  
  Observe that
  \begin{equation}
    \left| I_1 - 2 | \ln \eps | \, \| \msq \cdot \nnn\|^2_{L^2
    \left( \partial \Om \right)} \right| \leqslant \left| J_1 - 2 | \ln
    \eps | \, \| \msq \cdot \nnn \|^2_{L^2 \left( \partial \Om
    \right)} \right| + | J_2 | + | J_3 |, \label{eq:estI101}
  \end{equation}
  and we first want to estimate $J_2$ and $J_3$. Using the estimate in
  Lemma~\ref{Lemma:LemmaforK1}, we obtain that as soon as $\eps$ is
  small enough, there holds
  \begin{align}
    |J_2|
    & \leqslant   \int_{\Od} \int_{\Od} \frac{\left| \nabla \eta_{\eps}
                  (\rr) \right|  \left| \nabla \eta_{\eps} (\sss) \right|  \left| \msq
                  (\rr) - \msq (\bpi (\rr)) \right|}{| \rr - \sss |}  \hspace{0.17em}
                  \hspace{0.17em} \dd \rr  \dd \sss \nonumber\\
    & \eqs   \int_{\Od} f_{\eps} (\rr)  \left| \nabla \eta_{\eps} (\rr)
             \right|  \left| \msq (\rr) - \msq (\bpi (\rr)) \right| \hspace{0.17em} \dd
             \rr  \nonumber\\
    & \leqslant  3  | \ln \eps  |\left( \int_{\Od} \left| \nabla \eta_{\eps} (\rr)
                  \right| | \msq (\rr) - \msq (\bpi (\rr)) | \hspace{0.17em} \dd \rr \right)
                  .
  \end{align}
  Applying the $L^1$-type estimate in Lemma~\ref{prop:boundsforj2j3}
  with $\mu = \eps$, we infer that
  \begin{equation}
    |J_2| \leqslant C | \ln \eps | \left( \eps + | \ln \eps |
      \, \|
      \eta_{\eps} \nabla \tmmathbf{m}_{\bot}^{\eps} \|_{L^1 \left( \Od
        \right)} \right),
  \end{equation}
  for some $C > 0$. Using the Cauchy-Schwarz and Young's inequalities,
  we then obtain
  \begin{equation}
    |J_2|  \leqslant  C | \ln \eps |  \left( \eps + | \ln
      \eps | \, |\Od|^{1/2} \| \eta_{\eps} \nabla
      \tmmathbf{m}_{\bot}^{\eps} \|_{L^2 
        ( \Om_{\eps}),} \right) \leqslant C' \eps^{1 / 2} | \ln
    \eps |^2 \left( 1 + \| \eta_{\eps} \nabla \tmmathbf{m}_{\bot}^{\eps}
      \|^2_{L^2 ( \Om_{\eps})} \right) , \label{esJ2}
  \end{equation}
  for some $C' > 0$ and all $\eps$ small enough.  In the same way, we
  obtain
  \begin{equation}
    |J_3| \leqslant 2 C' \eps^{1 / 2} | \ln \eps |^2 \left( 1 + \|
      \eta_{\eps} \nabla \tmmathbf{m}_{\bot}^{\eps} \|^2_{L^2 \left(
          \Om_{\eps} \right)} \right) ,\label{esJ3}
  \end{equation}
  for all $\eps$ small enough.  Hence, from \eqref{eq:estI101},
  \eqref{esJ2} and \eqref{esJ3} we get that
  \begin{multline}
    \left| I_1 - 2 | \ln \eps | \| \msq \cdot \nnn \|^2_{L^2
        \left( \partial \Om \right)} \right| \\
    \leqslant \left| J_1 - 2 | \ln
    \eps | \, \| \msq \cdot \nnn \|^2_{L^2 \left( \partial \Om
    \right)} \right| + C \eps^{1 / 2} | \ln \eps |^2 \left( 1 + \|
    \eta_{\eps} \nabla \tmmathbf{m}_{\bot}^{\eps} \|^2_{L^2 \left(
    \Om_{\eps} \right)} \right), \label{eq:estI102} 
\end{multline}
for all $\eps$ small enough.

It remains to estimate
$|\, J_1 - 2 | \ln \eps | \, \| \msq \cdot \nnn \|^2_{L^2 \left(
      \partial \Om \right)} \, |$. We proceed by decomposing $J_1$
as $J_1 \assign K_1 + K_2$ with
  \begin{align}
    K_1 & \assign  \int_{\Od} \int_{\Od} \frac{\hspace{0.17em} \left| \nnn
                    (\bpi (\rr)) \cdot \msq (\bpi (\rr)) \right|^2}{| \rr - \sss |}  \left|
                    \nabla \eta_{\eps} (\rr) \right|  \left| \nabla \eta_{\eps} (\sss)
                    \right|  \hspace{0.17em} \dd \rr \dd \sss, \\
        &    \nonumber\\
    K_2 & \assign  \int_{\Od} \left| \nabla \eta_{\eps} (\rr) \right| \nnn
                    (\bpi (\rr)) \cdot \msq (\bpi (\rr))_{_{_{_{}}}} \nonumber\\
        &   \qquad \times \int_{\Od} \left| \nabla \eta_{\eps} (\sss)
             \right|  \frac{\nnn \left( \bpi (\sss) \right) \cdot \msq (\bpi (\sss)) -
             \nnn (\bpi (\rr)) \cdot \msq (\bpi (\rr))}{| \rr - \sss |} \dd \sss
             \dd \rr ,
  \end{align}
  and we show that
  \begin{eqnarray}
    \left| J_1 - 2 | \ln \eps | \left\| \msq \cdot \nnn \right\|^2_{L^2
    \left( \partial \Om \right)} \right|
    & \leqslant & \left| K_1 - 2 | \ln
                  \eps | \left\| \msq \cdot \nnn \right\|^2_{L^2
                  \left( \partial \Om 
                  \right)} \right| + | K_2 | \nonumber\\
    & \leqslant & C \left( 1 + \| \eta_\eps \nabla \msq
                  \|_{L^2(\Om)} \right) 
                  \left\| \msq \cdot \nnn \right\|_{L^2 \left( 
                  \partial \Om \right)} ,
                  \label{J1}
  \end{eqnarray}
  for some $C > 0$ and all $\eps$ small enough.
  
  We estimate $K_2$ to obtain
  \begin{multline}
    | K_2 |\leqslant \int_{\Od} \int_{\Od} \frac{\left| \msq (\bpi
        (\sss)) \cdot \nnn \left( \bpi (\sss) \right) - \msq (\bpi
        (\rr)) \cdot \nnn
        (\bpi (\rr)) \right|}{| \bpi (\rr) - \bpi (\sss) |}  \\
    \times \left| \nnn (\bpi (\rr)) \cdot \msq (\bpi (\rr)) \right|
    \frac{| \bpi (\rr) - \bpi (\sss) |}{| \rr - \sss |} \left| \nabla
      \eta_{\eps} (\rr) \right| \left| \nabla \eta_{\eps} (\sss)
    \right| \hspace{0.17em} \dd \rr \dd \sss .
  \end{multline}
  Since $\partial \Om$ is of class $C^2$ and compact, the projection
  map $\bpi : \Odbar \rightarrow \partial \Om$ is uniformly Lipschitz
  for sufficiently small $\bar \eps$. Thus, there exists a constant
  ${C_{\bpi}} > 0$ such that
  \begin{equation}
    | \bpi (\rr) - \bpi (\sss) | \leqslant C_{\bpi}  | \rr - \sss | \quad \forall
    \rr, \sss \in \Odbar ,
  \end{equation}
  and passing to the curvilinear coordinates we obtain
  \begin{align}
    | K_2 |
    & \leqslant  C_{\bpi} \int_{\Od} \int_{\Od} \frac{\left| \msq (\bpi
                  (\sss)) \cdot \nnn \left( \bpi (\sss) \right) - \msq (\bpi (\rr)) \cdot \nnn
                  (\bpi (\rr)) \right|}{| \bpi (\rr) - \bpi (\sss) |} \cdot \nonumber\\
    &   \qquad \qquad \qquad \cdot \left| \nnn (\bpi (\rr))
         \cdot \msq (\bpi (\rr)) \right| \left| \nabla \eta_{\eps} (\rr) \right| 
         \left| \nabla \eta_{\eps} (\sss) \right| \hspace{0.17em} \dd \rr  \dd
         \sss \\
    & \leqslant  2 C_{\bpi} \int_0^1 \int_0^1 | \eta' (s) | | \eta' (t) |
                  \int_{\partial \Om} \int_{\partial \Om} \frac{\left| \msq (\mu) \cdot \nnn
                  (\mu) - \msq (\sigma) \cdot \nnn (\sigma) \right|}{| \mu - \sigma |} \cdot
                  \nonumber\\
    &   \qquad \qquad \qquad \cdot \left| \nnn (\sigma) \cdot
         \msq (\sigma) \right|   \dd \mathcal H^1 (\mu)
         \dd \mathcal H^1(\sigma) \dd s \dd t, 
  \end{align}
  provided $\bar \eps$ is small enough. 
  
  Since $\| \eta' \|_{L^1 (0, 1)} = 1$, using the Cauchy-Schwarz
  inequality we obtain
  \begin{align}
    | K_2 |
    & \leqslant  C \int_{\partial \Om} \int_{\partial \Om}
                  \frac{\left| \msq (\mu) \cdot \nnn (\mu) - \msq (\sigma) \cdot \nnn
                  (\sigma) \right|}{| \mu - \sigma |}  \left| \nnn (\sigma) \cdot \msq
                  (\sigma) \right|  \dd \mathcal H^1(\mu) \dd \mathcal
                  H^1(\sigma) \; \nonumber\\ 
    & \leqslant  C' \left( \int_{\partial \Om} \left| \nnn (\sigma) \cdot
                  \msq (\sigma) \right|^2 \dd \mathcal H^1(\sigma)
                  \right)^{\frac{1}{2}}  \notag \\
    &\qquad \qquad \times \left(
       \int_{\partial \Om} \int_{\partial \Om} \frac{\left| \msq (\mu) \cdot \nnn
       (\mu) - \msq (\sigma) \cdot \nnn (\sigma) \right|^2}{| \mu - \sigma |^2} 
       \dd \mathcal H^1(\mu) \dd \mathcal H^1(\sigma)
       \right)^{\frac{1}{2}} \nonumber\\ 
    & \leqslant  C'' \| \msq \cdot \nnn \|_{H^{1 / 2} \left(
                  \partial \Om \right)} \| \msq \cdot \nnn \|_{L^2 \left(
                  \partial \Om \right)} ,  \label{esK2}
  \end{align}
  for some $C, C', C'' > 0$ and all $\eps$ small enough. Finally,
  using $\left| \msq \right| \leqslant 1$ and the trace embedding
  theorem, we obtain
  \begin{equation}
    | K_2 | \leqslant C \left( 1 +  \| \eta_{\eps}
      \nabla \msq \|_{L^2 ( \Om_{\eps})}  \right) 
    \| \msq \cdot \nnn \|_{L^2 ( \partial \Om )} .
    \label{estK2n}
  \end{equation}

  Lastly, we show that
  \begin{equation}
    \left| K_1 - 2 | \ln \eps | \int_{\partial \Om} ( \mpe^{\eps}
      (\sigma) \cdot \nnn (\sigma))^2  \dd \mathcal H^1(\sigma) \right| \leqslant C
    \| \msq \cdot \nnn\|_{L^2( \partial \Om)} .
    \label{eqK1}
  \end{equation}
  Indeed, using the coarea formula and recalling the definition of
  $f_\eps$ in \eqref{eq:fdeltanew0}, we have
  \begin{equation}
    K_1 
     \eqs \int_{\partial \Om} \left| \nnn (\sigma) \cdot \msq (\sigma)
             \right|^2 \left( \int_0^1 | \eta' (t) | f_{\eps} (\sigma + \eps
             t\tmmathbf{n} (\sigma))  ( 1 + \eps t \kappa (\sigma) )  \hspace{0.17em}
             \dd t \right)  \dd \mathcal H^1(\sigma) . 
  \end{equation}
  Therefore, using the asymptotics of
  $f_{\eps} (\sigma + \eps t\tmmathbf{n} (\sigma))$ given in
  Lemma~\ref{Lemma:LemmaforK1} and the fact that
  $\left| \msq \right| \leqslant 1$, we infer \eqref{eqK1}. Combining
  \eqref{estK2n} and \eqref{eqK1}, we get \eqref{J1}.
  Finally, combining it with \eqref{eq:estI102}, \eqref{eq:firstsubd}
  \eqref{estI3n}, and \eqref{eq:VOmcrossOmplusI2} we get the desired
  estimate \eqref{eq:2beproved}. This concludes the proof.
\end{proof}

\begin{proposition}
  \label{th:mainresultdemag2} There exist $\dbar, C > 0$ such that if
  $0 < \eps < \dbar$ and $\m^\eps \in \myspace$, then the
  magnetostatic energy for the out-of-plane component
  {\tmem{(}}{\tmabbr{cf.}}~{\tmem{\eqref{eq:exprmag2})}}
  \begin{align}
    \tilde{\mathcal{V}} ( \eta_{\eps} \mpa^{\eps}) =
    \int_{\R^2} \int_{\R^2} \frac{\nabla \left( \eta_{\eps}
    \mpa^{\eps} \right) ( \rr ) \cdot \nabla \left( \eta_{\eps}
    \mpa^{\eps} \right) ( \sss) }{| \rr - \sss | } \hspace{0.17em}
    \dd \rr \dd \sss
  \end{align}
  satisfies
   \begin{align}
     \def\arraystretch{1.4}
     \begin{array}{l} \left|
       \tilde{\mathcal{V}}( \eta_{\eps} \mpa^{\eps} ) -
       \tilde{\mathcal{V}}_{\Om \times \Om} (\mps) + 2
       \tilde{\mathcal{V}}_{\partial \Om \times \Om} ( \mps ) -
       \Ddelta +\displaystyle 2 | \ln \eps | \int_{\partial \Om}
       \left( 1 - | \mpa^{\eps}|^2 \right) \dd \mathcal H^1(\sigma) \right|
       \\
       \qquad \qquad \leqslant 
       C_{\eps} \left( 1 + \|
       \eta_{\eps} \nabla \mps
       \|_{L^2 ( \Om_{\eps}
       )}^2 \right) + C \|
       \mps + 1\|^{1 / 2}_{L^2 (
       \partial \Om )} \| \mps - 1
       \|^{1 / 2}_{L^2 \left(
       \partial \Om \right)} \left(
       1 + \| \eta_{\eps} \nabla
       \mps \|_{L^2 ( \Om_{\eps} )}
       \right) , 
     \end{array}
   \end{align}
  where $C_{\eps} \rightarrow 0$ as $\eps \rightarrow 0$, $D_\eps$
    is defined in \eqref{eq:Ddeltadef}, and
  \begin{align}
    \tilde{\mathcal{V}}_{\Om \times \Om}( \mps ) & \assign 
    \int_{\Om} \int_{\Om}  \frac{\nabla \mps (\rr) \cdot
    \nabla \mps ( \sss)}{\left| \sss - \rr \right|}  \dd \sss  \dd
    \rr, \\
    \tilde{\mathcal{V}}_{\partial \Om \times \Om} ( \mps ) &
    \assign \int_{\Om} \int_{\partial \Om}  \frac{\nabla
    \mps (\rr) \cdot \nnn (\sigma) \mps (\sigma)}{\left| \sigma - \rr
                                                             \right|} 
    \dd \mathcal H^1(\sigma) \dd \rr . 
  \end{align}
\end{proposition}

\begin{proof}
  We begin by writing $\tilde{\mathcal{V}}$ in the form similar to
  that of the nonlocal term in Proposition \ref{th:mainresultdemag}:
  \begin{align} \tilde{\mathcal{V}} ( \eta_{\eps} \mpa^{\eps}) =
    \sum_{i = 1}^2 \int_{\R^2} \int_{\R^2} \frac{\divv \left(
        \eta_{\eps} \mpa^{\eps} \tmmathbf{e}_i \right) ( \rr ) \divv
      \left( \eta_{\eps} \mpa^{\eps} \tmmathbf{e}_i \right) ( \sss)
    }{| \rr - \sss | } \hspace{0.17em} \dd \rr \dd \sss . \end{align}
  Proceeding exactly as in {\tmname{Steps}}~1--3 in the proof of
  Proposition~\ref{th:mainresultdemag}, we obtain
  \begin{align} \left| \tilde{\mathcal{V}} ( \eta_{\eps} \mpa^{\eps})
      - \tilde{J}_1 - \tilde{\mathcal{V}}_{\Om \times \Om} ( \mps) + 2
      \tilde{\mathcal{V}}_{\partial \Om \times \Om} ( \mps ) \right|
    \leqslant C_\eps ( 1 + \| \eta_{\eps} \nabla \mps \|^2_{L^2 (
      \Om_{\eps})} ),
  \end{align}
  where $C_\eps \to 0$ as $\eps \to 0$ and
  \begin{align} 
       \tilde{J}_1 \assign \int_{\Od} \int_{\Od} \frac{\nabla
       \eta_{\eps} (\rr) \cdot \nabla \eta_{\eps} (\sss)  \mps (\bpi (\rr))
       \hspace{0.17em} \mps (\bpi (\sss))}{| \rr - \sss |}  \hspace{0.17em} \dd
       \rr  \dd \sss .
  \end{align}
  
  To account for the non-zero limiting boundary data for $\hspace{0.17em}
  \mps$, we represent $\tilde{J}_1$ in the following way:
  \begin{align}
    \tilde{J}_1 = \int_{\Od} \int_{\Od} \frac{\nabla \eta_{\eps} (\rr)
    \cdot \nabla \eta_{\eps} (\sss)  \left( \mps (\bpi (\rr)) - 1 \right) 
    \hspace{0.17em} \left( \mps (\bpi (\sss)) + 1 \right)}{| \rr - \sss |} 
    \hspace{0.17em} \dd \rr  \dd \sss + \Ddelta
  \end{align}
  where we recalled the definition of $D_\eps$ from
  \eqref{eq:Ddeltadef} and noted that the terms linear in $\mps$
  cancel upon expansion of the integral due to the symmetry of the
  kernel. We denote by $\bar{J}_1$ the first integral in the above
  expression (i.e., $\bar{J}_1 \assign \tilde{J}_1 - \Ddelta$) and
  split it in a similar way to what we did for $J_1$ in
  Proposition~\ref{th:mainresultdemag}. Specifically, we set
  $\bar{J}_1 \assign \tilde{K}_1 + \tilde{K}_2$, where
  \begin{align}
    \tilde{K}_1 & \assign  \int_{\Od} \int_{\Od} \frac{\hspace{0.17em} \left|
    \mps (\bpi (x)) \right|^2 - 1}{| \rr - \sss |}  \left| \nabla \eta_{\eps}
    (\rr) \right|  \left| \nabla \eta_{\eps} (\sss) \right|  \hspace{0.17em}
    \dd \rr \dd \sss, \\
    \tilde{K}_2 & \assign  \int_{\Od} \left| \nabla \eta_{\eps} (\rr)
    \right| \left( \mps (\bpi (x)) - 1 \right) \nnn (\bpi (x))_{_{_{_{}}}}
    \nonumber \\
    &   \qquad \cdot \int_{\Od} \left| \nabla \eta_{\eps} (\sss) \right| 
    \frac{\nnn \left( \bpi (\sss) \right) \left( \mps (\bpi (y)) + 1 \right) -
    \nnn (\bpi (\rr)) \left( \mps (\bpi (\rr)) + 1 \right)}{| \rr - \sss |}
    \dd \sss \dd \rr . 
  \end{align}
  By the same arguments used in the proof of {\tmname{Step~3}} in
  Proposition~\ref{th:mainresultdemag} to estimate $K_2$, we then
  obtain the estimate
  \begin{equation}
    \tilde{K}_2 \leqslant C \left( 1 + \| \eta_{\eps} \nabla \mps
      \| _{L^2 ( \Om_{\eps})} \right) \| \mps - 1
    \|_{L^2 ( \partial \Om )},   \label{K2parx}
\end{equation}
for some $C > 0$ and all $\eps$ small enough.

Alternatively, writing $\tilde{K}_2$ in the following equivalent way:
  \begin{align}
    \tilde{K}_2 & \eqs  \int_{\Od} \left| \nabla \eta_{\eps} (\rr) \right|
    \left( \mps (\bpi (x)) + 1 \right) \nnn (\bpi (x))_{_{_{_{}}}} \nonumber\\
    &  \qquad \cdot \int_{\Od} \left| \nabla \eta_{\eps} (\sss) \right| 
    \frac{\nnn \left( \bpi (\sss) \right) \left( \mps (\bpi (y)) - 1 \right) -
    \nnn (\bpi (\rr)) \left( \mps (\bpi (\rr)) - 1 \right)}{| \rr - \sss |}
    \dd \sss \dd \rr , 
  \end{align}
  we infer
  \begin{equation}
    \tilde{K}_2 \leqslant C
    \left( 1 + \| \eta_{\eps} \nabla \mps
      \| _{L^2 ( \Om_{\eps})} \right) \| \mps + 1 \|_{L^2 (
      \partial \Om)},  \label{K2pary}
  \end{equation}
  and taking the geometric mean of \eqref{K2parx} and \eqref{K2pary},
  we obtain
  \begin{equation}
    \tilde{K}_2 \leqslant C \left( 1 + \| \eta_{\eps} \nabla \mps
      \| _{L^2 ( \Om_{\eps})} \right)\| \mps + 1
    \|^{1 / 2}_{L^2 ( \partial \Om )} \, \| \mps - 1
    \|^{1 / 2}_{L^2 ( \partial \Om)} , \label{eq:tK2n}
  \end{equation}
  for some $C > 0$ and all $\eps$ small enough.

  Finally, the estimate for $\tilde K_1$ can be obtained in the the
  same way we derived the estimate for $K_1$ in {\tmname{Step~3}} of
  the proof of Proposition~\ref{th:mainresultdemag}, together with a
  Cauchy-Schwarz inequality, to obtain
  \begin{equation}
    \left| \tilde{K}_1 - 2 | \ln \eps | \int_{\partial \Om}  ( |
      \mps|^2 - 1 )  \dd \mathcal H^1(\sigma) \right| \leqslant C \, \| \mps +
    1 \|_{L^2 ( \partial \Om )}^{1/2} \, \| \mps - 1
    \|_{L^2( \partial \Om)}^{1/2} . \label{tildek1}
  \end{equation}
  Combining all of the above estimates, we obtain the result.
\end{proof}

\section{Proof of
  \texorpdfstring{$\Gamma$}{Gamma}-convergence}\label{sec:Gammaconv}

In this section we provide the proof of our main theorems. Since the
proofs of Theorem \ref{t:GJ} and Theorem~\ref{thm:locengy} follow
essentially verbatim those of Theorem~\ref{thm:thmKS} and
Theorem~\ref{thm:nonlocengy}, respectively, we only give the proofs of
the latter. Theorem~\ref{t:GJ} and Theorem~\ref{thm:locengy} may in
fact be thought of as the limiting cases of Theorem~\ref{thm:thmKS}
and Theorem~\ref{thm:nonlocengy}.

\begin{proof}[Proof of Theorem~\ref{thm:thmKS}] Without loss
  of generality we may suppose that $\gamma_\eps = \gamma$. \\
  {\tmem{}}(i) (Compactness) We first prove the compactness
  result. Let us assume that
  $\mathcal{G}_{\eps} (\m^{\eps}) \leqslant C$ for some constant
  $C > 0$ independent of $\eps$. We recall that
  ({\tmabbr{cf.}}~\eqref{eq:enGdelta})
  \begin{equation}
    \mathcal{G}_{\eps} (\m^{\eps}) \eqs  \|
    \eta_{\eps} \nabla \m^{\eps} \|_{L^2 \left( \Om_{\eps}
    \right)}^2 + \lambda \DMIdelta (\m^{\eps}) + \frac{\gamma}{2 | \ln
    \eps |} \Vdelta ( \eta_{\eps} \msq) - \frac{\gamma}{2 |
    \ln \eps |} \Vtildedelta ( \eta_{\eps} \mpa^{\eps}) .
  \end{equation}
  First, note that, up to a constant term, we can absorb the DMI
  energy $\lambda \DMIdelta (\m^{\eps})$ into the Dirichlet
  energy. Indeed, since $\left| \m \right| = 1$ a.e. in $\Om_{\eps}$
  and $| \eta_{\eps} | \leqslant 1$, by the Cauchy-Schwarz and Young's
  inequalities for every $0 < \varepsilon < \dbar$ and every
    $\delta > 0$ there holds
  \begin{align}
    \left| \DMIdelta (\m^{\eps}) \right|
    & \leqslant
      \int_{\Om_{\eps}} \left| \mpa^{\eps} \divv \mpr^{\eps} -
      \mpr^{\eps} \cdot \nabla \mpa^{\eps} \right| | \eta_{\eps} | \dd
      \rr \notag \\
    & \leqslant C \| \eta_{\eps} \nabla \m^{\eps} \|_{L^2(
      \Om_{\eps} )} \leqslant \frac{\delta}{2} \| \eta_{\eps}
      \nabla \m^{\eps} \|_{L^2( \Om_{\eps})}^2 + \frac{C^2}{2
      \delta},
  \end{align}
  for some $\dbar, C > 0$ that depend only on $\Om$.
  Therefore, without loss of generality, we can assume from the very
  beginning that
  \begin{align}
    \label{etamVVC}
    \| \eta_{\eps} \nabla \m^{\eps}\|_{L^2(
    \Om_{\eps})}^2 + \frac{\gamma}{ | \ln \eps |} \Vdelta(
    \eta_{\eps} \msq) - \frac{\gamma}{ | \ln \eps |} \Vtildedelta
    ( \eta_{\eps} \mpa^{\eps}) \leqslant C
  \end{align}
  for some constant $C > 0$ independent of $\eps$. 
  
  By positivity of $\Vdelta( \eta_{\eps} \msq)$, we may, furthermore,
  drop this term from \eqref{etamVVC}.  On the other hand, from
  Proposition \ref{th:mainresultdemag2}, Lemma \ref{lem:Deps} and
    the estimates
    \begin{align}
      \tilde{\mathcal V}_{\Om \times \Om} (\mpa^{\eps})
      & \leqslant \int_\Om \int_\Om \frac{|\nabla \mpa^\eps(x)|^2}{ 
        |x - y|} \, \dd y \, \dd x \leqslant C' \| \eta_{\eps}
        \nabla \m^{\eps}\|_{L^2( \Om_{\eps})}^2 , \label{VOmOm} \\
           \tilde{\mathcal V}_{\partial \Om \times \Om} (\mpa^{\eps})
      & \leqslant \int_\Om \int_{\partial \Om} \frac{|\nabla \mpa^\eps(x)}{
        |x - \sigma|} \, \dd \mathcal H^1(\sigma) \, \dd x
        \leqslant C'' \| \eta_{\eps} 
        \nabla \m^{\eps}\|_{L^2( \Om_{\eps})} , \label{VdOmOm} 
    \end{align}
    for some $C', C'' > 0$ depending only on $\Omega$ that follow from
    the Cauchy-Schwarz inequalities, we immediately obtain the
    existence of a positive constant $C > 0$ such that for all
    sufficiently small $\eps$ there holds
  \begin{equation}
    (1 - \gamma C_{\eps}) \, \| \eta_{\eps}
      \nabla \m^{\eps}\|_{L^2( \Om_{\eps})}^2 \leqslant C
    , \label{1gCepsC} 
  \end{equation}
  for some $C_{\eps} > 0$ such that $C_{\eps} \rightarrow 0$ as
  $\eps \rightarrow 0$.
  
  From \eqref{1gCepsC} we conclude that for $\eps$ sufficiently small
  we can uniformly bound \ $\m^{\eps}$ in
  $H^1 \left( \Om; \Stwo^2 \right)$ and, therefore, up to a
  subsequence, there exists $\m \in H^1 \left( \Om; \Stwo^2 \right)$
  such that
  \begin{align}
    \m^{\eps}  & \rightarrow  \m
                 \quad \text{$\text{strongly in }
                 L^2(\Om; \RR^3 )$} \\
    \nabla \m^{\eps}  & \rightharpoonup \nabla \m
                        \quad \text{weakly in }
                        L^2 ( \Om; \RR^{2 \times 3}), \\
    \mps  & \rightarrow  \mpa
            \quad \text{strongly in $L^p(
            \partial \Om )$ for $p \geqslant 1$},  \label{mpabdcon}\\
    \msq  & \rightarrow  \mpe
            \quad \text{strongly
            in $L^p \left( \partial \Om; \RR^2 \right)$ for $p
            \geqslant 1$} , 
    \label{plbdcon}
  \end{align}
  as $\eps \to 0$.  Moreover, since $\left| \m^{\eps} \right| = 1$ on
  $\Om_{\eps}$ we also have $\eta_{\eps} \m^{\eps} \rightarrow 0$
  strongly in $L^2 ( \R^2 \setminus \overline{\Om}; \RR^3 )$. Hence,
  we infer that if $\m$ is extended by zero outside $\Omega$ then
  \begin{equation}
    \m^{\eps} \rightarrow \m \quad \text{strongly
      in } L^2 ( \RR^2; \RR^3) .
  \end{equation}
  In particular, we have $\m \in H_0(\RR^2; \mathbb
  S^2)$. \\ \\
  {\tmem{}}(ii) ({\tmem{$\Gamma \textrm{\text{-}liminf}$}} inequality)
  Let $\m^\eps \in \myspace$ and $\m \in H_0(\RR^2; \mathbb S^2)$ be
  such that $\m^{\eps} \rightarrow \m$ strongly in
  $L^2 ( \RR^2; \RR^3)$ as $\eps \rightarrow 0$. We may further assume
  that
  $\liminf_{\eps \rightarrow 0} \mathcal{G}_{\eps} (\m^{\eps}) < +
  \infty$, since otherwise the statement is trivially true. Hence,
  using the compactness statement (maybe passing to a subsequence) we
  have $\m_{\eps} \rightharpoonup \m$ weakly in
  $H^1 \left( \Om; \RR^3 \right)$, and using
  Propositions~\ref{th:mainresultdemag} and \ref{th:mainresultdemag2},
  together with Lemma \ref{lem:Deps} and the lower semicontinuity of
  the Dirichlet energy on $\Om$ and the compactness of trace embedding
  of functions in $H^1( \Om)$ into $L^2 \left( \partial \Om \right)$
  we obtain
  \begin{multline}
    \liminf_{\eps \rightarrow 0} \mathcal{G}_{\eps} \left( \m^{\eps} 
    \right)
    \geqslant  \int_{\Om} \left| \grad \m \right|^2  \dd \rr +
      \lambda \int_{\Om} \left( \mpa  \divv  \mpr - \mpr \cdot \nabla
      \mpa \right)  \dd \rr \\
    + \gamma \int_{\partial \Om} \left( (\mpr \cdot
      \nnn)^2 - \mpa^2 \right)  \dd \mathcal H^1(\sigma)
      \eqs \mathcal{G}_0 \left( \m  \right) . 
    \end{multline}
    {\noindent}{\smallskip}(iii)
    ({\tmem{$\Gamma \textrm{\text{-}limsup}$}} inequality) Let
    $\m \in H_0 \left( \RR^2; \Stwo^2 \right)$ be such that
    $\mathcal{G}_0 \left( \m \right) < + \infty$. Take
    $\bar{\eps} > 0$ sufficiently small and extend $\m$ to
    $\widetilde{\m} \in H^1( \Om_{\bar{\eps}}, \Stwo^2)$, e.g., by
    setting
    $\widetilde{\m}(x) := \m(x - 2 d_{\partial \Om}(x)
    \nnn(\bpi(x)))$. For every $\eps < \bar{\eps}$ we now define
    $\m^{\eps} = \widetilde{\m}$ in $\Om_\eps$ and $\m^\eps = 0$
      outside $\Om_\eps$. It is clear that $\m^{\eps} \in \myspace$
    and $\m^\eps \rightarrow \m$ strongly in
    $L^2( \RR^2; \RR^3)$ as $\eps \rightarrow 0$. Moreover, using
    Propositions~\ref{th:mainresultdemag} and
    \ref{th:mainresultdemag2}, Lemma \ref{lem:Deps} and the
    strong convergence of $\m^{\eps} $ to $\m$ in
    $H^1( \Om; \RR^3)$, we have $\m^{\eps} \rightarrow \m$ in
    $L^2( \partial \Om; \RR^3)$ and can pass to the limit in the
    magnetostatic energy term. Finally, using the fact that
  \begin{align} \int_{\Om_{\eps}} \eta_{\eps }^2 \left| \nabla \m^{\eps} \right|^2 
     \dd \rr = \int_{\Od} \eta_{\eps }^2  \left| \nabla \widetilde{\m}
     \right|^2 \dd \rr + \int_{\Om} \left| \nabla \m \right|^2  \dd \rr \; \;
     \xrightarrow{\eps \rightarrow 0} \; \; \int_{\Om} \left| \nabla \m
     \right|^2  \dd \rr, \end{align}
  we obtain
  \begin{align} \limsup_{\eps \rightarrow 0} \mathcal{G}_{\eps} \left( \m^{\eps}
     \right) =\mathcal{G}_0 \left( \m  \right) . \end{align}
  This completes the proof.
\end{proof}

\begin{proof}[Proof of Theorem~\ref{thm:nonlocengy}]
  (i) (Compactness) We first prove the compactness result. Let us
  assume that
  $\mathcal{G}_{\eps} (\m^{\eps}) + \frac{\nu}{2} \Ddelta \leqslant C$
  for some constant $C > 0$ independent of $\eps$. We recall that now
  $\mathcal{G}_{\eps}$ reads as ({\tmabbr{cf.}}~\eqref{eq:enGdelta}
  with $\gamma_\eps = \nu | \ln \eps |$)
  \begin{equation}
    \mathcal{G}_{\eps} (\m^{\eps})  \eqs  \|
    \eta_{\eps} \nabla \m^{\eps} \|_{L^2( \Om_{\eps}
    )}^2 + \lambda \DMIdelta (\m^{\eps}) + \frac{\nu}{2} \Vdelta
    ( \eta_{\eps} \msq) - \frac{\nu}{2} \Vtildedelta \left(
    \eta_{\eps} \mpa^{\eps} \right) . 
  \end{equation}
  As in the proof of Theorem~\ref{thm:thmKS}, up to a constant term,
  we can absorb the DMI energy $\lambda \DMIdelta (\m^{\eps})$ into
  the Dirichlet energy and drop the $\Vdelta ( \eta_{\eps} \msq)$
    term due to its positivity. Therefore, without loss of
  generality, we can assume from the very beginning that
  \begin{align}
    \| \eta_{\eps} \nabla \m^{\eps} \|_{L^2(
     \Om_{\eps} )}^2 - \nu \Vtildedelta ( \eta_{\eps} \mpa^{\eps}
    ) + \nu \Ddelta \leqslant C, \label{eq:estimthm5proof1}
  \end{align}
  for some constant $C > 0$ independent of $\eps$.

  Next, aiming to invoke Proposition~\ref{th:mainresultdemag2}, we
  write
  \begin{align}
    - \nu (\Vtildedelta(
    & \eta_{\eps} \mpa^{\eps}    )
      - \Ddelta )
      \eqs  \nu \left( 2
      \tilde{\mathcal{V}}_{\partial \Om \times \Om} ( \mps ) -
      \tilde{\mathcal{V}}_{\Om \times \Om} ( \mps ) + 2 | \ln \eps
      | \int_{\partial \Om} \left( 1 - | \mpa^{\eps} |^2 \right) \dd \mathcal H^1(\sigma)
      \right) \nonumber\\
    &   - \nu \left( \Vtildedelta ( \eta_{\eps}
      \mpa^{\eps} ) - \Ddelta - \tilde{\mathcal{V}}_{\Om \times \Om}
      ( \mps ) + 2 \tilde{\mathcal{V}}_{\partial \Om \times \Om}
      ( \mps) + 2 | \ln \eps | \int_{\partial \Om} \left( 1 - |
      \mpa^{\eps} |^2 \right) \dd \mathcal H^1(\sigma) \right),
  \end{align}
  from which by Proposition~\ref{th:mainresultdemag2} it follows
  that
  \begin{align}
    - \nu (\Vtildedelta( \eta_{\eps} \mpa^{\eps}
    ) - \Ddelta )  \geqslant
    &  -C - C' \| \eta_{\eps} \nabla \mps
      \|_{L^2( \Om_{\eps})} - C_{\eps} \| \eta_{\eps} \nabla \mps
      \|^2_{L^2( \Om_{\eps})} \notag \\
    & + 2 \nu | \ln \eps |
      \int_{\partial \Om} \left( 1 - | \mpa^{\eps} |^2 \right)  \dd
      \sigma 
      + \nu \left( 2 \tilde{\mathcal{V}}_{\partial \Om \times
      \Om} ( \mps ) - \tilde{\mathcal{V}}_{\Om \times \Om} (
      \mps ) \right) , \label{nuVDeps}
  \end{align}
  where $C_\eps \to 0$ as $\eps \to 0$ and $C, C' > 0$ are independent
  of $\eps$. Now, from \eqref{VdOmOm} and Young's inequality it is
  clear that for any $\delta > 0$ we have
  \begin{equation}
    \nu \left| \tilde{\mathcal{V}}_{\partial \Om \times \Om} ( \mps )
    \right| 
    \leqslant 
    \delta \| \eta_{\eps} \nabla \mps \|^2_{L^2(
      \Om_{\eps})} + C \nu^2 \delta^{-1},
  \end{equation}
  for some $C > 0$ depending only on $\Om$.  Also, applying
  Lemma~\ref{lm:con} to $\uuu_\eps = \mpa^\eps \tmmathbf e_i$,
  $i = 1,2$, one obtains that for any $\delta > 0$ there holds
  \begin{align}
    \nu \left| \tilde{\mathcal{V}}_{\Om \times \Om} \left( \mps
    \right) \right| 
    \leqslant \delta \left\| \eta_{\eps} \nabla \mps
    \right\|^2_{L^2 ( \Om_{\eps}),} + C \nu^2 \delta^{-1} ,
  \end{align}
  again, for some $C > 0$ depending only on $\Om$.  Based on the above
  estimates and another application of Young's inequality in
  \eqref{nuVDeps} we deduce that for any $0 < \eps < \bar{\eps}$ and
  any $\delta > 0$ there holds
  \begin{equation}
    - \nu ( \Vtildedelta ( \eta_{\eps} \mpa^{\eps}
    ) - \Ddelta )  \geqslant - (  4 \delta + C_{\eps})
    \| \eta_{\eps} \nabla \mps \|^2_{L^2 ( \Om_{\eps}
      )} - C \left( 1 + \nu^2 \delta^{-1} \right)  + 2 \nu |
    \ln \eps 
    | \int_{\partial \Om} \left( 1 - | \mps|^2 \right)  \dd \mathcal H^1(\sigma) ,
  \end{equation}
  for some $C > 0$ depending only on $\Om$.  Therefore, we can
  absorb the term
  $-\nu \left( \Vtildedelta ( \eta_{\eps} \mpa^{\eps} ) - \Ddelta
  \right)$ into the Dirichlet energy by choosing $\delta$
    sufficiently small universal, and for any $0 < \eps < \bar{\eps}$
  there holds
  \begin{equation}
    \frac12 \| \eta_{\eps} \nabla \mps \|^2_{L^2 ( \Om_{\eps}
      )} + 2 \nu | \ln \eps | \int_{\partial
      \Om} \left( 1 - \left| \mpa^{\eps} \right|^2 \right) \dd
    \mathcal H^1(\sigma)
    \leqslant C,
  \end{equation}
  for some $C > 0$ independent of $\eps$.  This gives us (as in the
  proof of Theorem~\ref{thm:thmKS}) the existence of
  $\m \in H_0( \RR^2; \Stwo^2 )$ and a subsequence such that
  \begin{align}
    \m^{\eps} \rightarrow \m
    & \quad \text{strongly in } L^2 (
      \R^2; \RR^3 ), \\
    \nabla \m^{\eps}   \rightharpoonup \nabla \m
    & \quad
      \text{weakly
      in 
      } L^2( \Om; \RR^{2 \times 3} ).  \label{wgconpar}
    \\
    | \mps|  \rightarrow  1 & \quad \text{strongly in $L^p
                              ( \partial \Om )$ for any $p \geqslant 1$},  \label{mpabdconpar}\\
    \msq  \rightarrow  0 & \quad \text{strongly in $L^p
                                  ( \partial \Om; \R^2 )$ for any $p \geqslant
                                  1$,}  \label{plbdconpar} 
  \end{align}
  Hence, upon a further subsequence, we have $|\mpa^\eps| \to 1$
  a.e. in $\partial \Om$. In fact, since the trace of the limit
  belongs to $\text{VMO}(\partial \Om)$ and takes only values $\pm 1$,
  it is in fact constant a.e. on $\partial \Om$
  \cite{brezis95,brezis02}.
  
  In what follows, without loss of generality, we assume that
  $\mpa^\eps \to 1$ strongly in $L^p( \partial \Om )$,
  $p \geqslant 1$, i.e., that the limit configuration $\m$ satisfies
  the boundary condition $\m =\tmmathbf{e}_3$ a.e. on
  $\partial \Om$. \\
  
  {\noindent}(ii) ({\tmem{$\Gamma \textrm{\text{-}liminf}$}} inequality) We
  consider the energy functional
  \begin{align} \mathcal{G}_{\eps} (\m^{\eps}) + \frac{\nu}{2}
    \Ddelta \eqs \| \eta_{\eps} \nabla \m^{\eps} \|_{L^2 ( \Om_{\eps}
      )}^2 + \lambda \DMIdelta (\m^{\eps}) + \frac{\nu}{2} \Ddelta +
    \frac{\nu}{2} \Vdelta ( \eta_{\eps} \msq) - \frac{\nu}{2}
    \Vtildedelta ( \eta_{\eps} \mpa^{\eps}) .
  \end{align}
  and will prove a $\liminf$ inequality for this functional. Let
  $\m^{\eps} \in \myspace$ satisfy $\m^{\eps} \rightarrow \m$ strongly
  in $L^2 ( \RR^2; \RR^3)$ as $\eps \rightarrow 0$. We may assume that
  $\liminf_{\eps \rightarrow 0} \left( \mathcal{G}_{\eps} \left(
      \m^{\eps} \right) + \frac{\nu}{2} \Ddelta \right) < + \infty$,
  otherwise the statement is trivially true. Hence (maybe after
  passing to a subsequence) we may assume that 
  \begin{equation}
    \liminf_{\eps \rightarrow 0} \mathcal{G}_{\eps} \left( \m^{\eps}
    \right) + \frac{\nu}{2} \Ddelta = \lim_{\eps \rightarrow 0}
    \mathcal{G}_{\eps} (\m^{\eps}) + \frac{\nu}{2} \Ddelta <
    + \infty .
  \end{equation}
  Using the compactness result and the implied convergence, by the
  lower semicontinuity of the norm and the weak-strong argument we
  immediately obtain
  \begin{equation}
    \liminf_{\eps \rightarrow 0} \| \eta_{\eps} \nabla \m^{\eps}
    \|_{L^2 ( \Om_{\eps})}^2 \geqslant\|
    \nabla \m \|_{L^2 ( \Om)}^2, \qquad
    \lim_{\eps \rightarrow 0} \lambda \DMIdelta (\m^{\eps}) \eqs \lambda
    \int_{\Om} \left( \mpa  \divv  \mpr - \mpr \cdot \nabla \mpa \right) \dd
    \rr .
  \end{equation}
  Therefore, the {\tmem{$\Gamma \textrm{\text{-}liminf}$}} inequality is proved
  once we show that
  \begin{equation}
    \liminf_{\eps \rightarrow 0} \left( \frac{\nu}{2} \Ddelta +
    \frac{\nu}{2} \Vdelta ( \eta_{\eps} \msq) - \frac{\nu}{2}
    \Vtildedelta ( \eta_{\eps} \mpa^{\eps}) \right) 
    \geqslant \frac{\nu}{2} \mathcal{V}_{\Om \times \Om} \left( \mpe \right)
    - \frac{\nu}{2} \tilde{\mathcal{V}}_{\Om \times \Om} \left( \mpa \right) +
    \nu \tilde{\mathcal{V}}_{\partial \Om \times \Om} \left( \mpa \right)\, .
     \label{eq:tocompletegammaliminf}
  \end{equation}
  For that, we consider separately the convergence of the terms
    due to the in-plane and the out-of-plane components.
 
  \vspace{8pt}
    {\noindent}{\tmname{The in-plane magnetostatic contribution}}. As
    a direct consequence of Proposition~\ref{th:mainresultdemag},
    Lemma~\ref{lm:con}, and the convergence in \eqref{wgconpar},
    \eqref{plbdconpar}, we obtain
  \begin{equation}
    \frac{\nu}{2} \mathcal{V} ( \eta_{\eps} \msq) - \nu | \ln
    \eps | \int_{\partial \Om}  \left( \msq \cdot \nnn \right)^2  \dd
    \mathcal H^1(\sigma)
    \; \; \; \xrightarrow{\eps \rightarrow 0} \; \; \; \frac{\nu}{2}
    \mathcal{V}_{\Om \times \Om} \left( \mpe \right) .
    \label{eq:GammliminfDdelta1}
  \end{equation}
 
 \vspace{8pt} 
  {\noindent}{\tmname{The out-of-plane magnetostatic
      contribution.}} As a direct consequence of
  Proposition~\ref{th:mainresultdemag2}, Lemma~\ref{lm:con}, and the
  convergence in \eqref{wgconpar}, \eqref{mpabdcon}, we obtain
  \begin{equation}
    \frac{\nu}{2} \left( \Vtildedelta \left( \eta_{\eps} \mpa^{\eps}
      \right) - \Ddelta \right) + \nu | \ln \eps | \int_{\partial \Om}
    \left( 1 -
      | \mps |^2 \right) \dd \mathcal H^1(\sigma) \; \; \;
    \xrightarrow{\eps \rightarrow 
      0} \; \; \frac{\nu}{2} \tilde{\mathcal{V}}_{\Om \times \Om} ( \mpa
    ) - \nu \tilde{\mathcal{V}}_{\partial \Om \times \Om} ( \mpa
    ) . \label{eq:GammliminfDdelta2}
  \end{equation}
  {\tmname{{\noindent}Proof of \eqref{eq:tocompletegammaliminf}.}}
  Combining \eqref{eq:GammliminfDdelta1} and
  \eqref{eq:GammliminfDdelta2} we get
  \begin{eqnarray}
    \frac{\nu}{2} \Ddelta + \frac{\nu}{2} \Vdelta ( \eta_{\eps} \msq
    ) - \frac{\nu}{2} \Vtildedelta ( \eta_{\eps} \mpa^{\eps}
    ) & \geqslant & \frac{\nu}{2} \Vdelta ( \eta_{\eps} \msq) - \nu |
               \ln \eps | \int_{\partial \Om}  ( \msq \cdot \nnn )^2
               \dd \mathcal H^1(\sigma)
               \nonumber\\
      &  & \qquad - \frac{\nu}{2} ( \Vtildedelta ( \eta_{\eps}
           \mpa^{\eps} ) - \Ddelta ) - \nu | \ln \eps |
           \int_{\partial \Om} \left( 1 -| \mps |^2 \right) \dd \mathcal
           H^1(\sigma) \notag \\ 
      & \xrightarrow{\eps \rightarrow 0} & \frac{\nu}{2} \mathcal{V}_{\Om
                                           \times \Om} ( \mpe ) - \frac{\nu}{2} \tilde{\mathcal{V}}_{\Om
                                           \times \Om} ( \mpa) + \nu \tilde{\mathcal{V}}_{\partial \Om
                                           \times \Om} ( \mpa ) .  \label{eq:ineqgammaliming}
  \end{eqnarray}
  Therefore, using the definition of the vector field
  $\tmmathbf{b}(x)$ (see \eqref{eq:b}) we obtain
  \begin{align} \liminf_{\eps \rightarrow 0} \left( \mathcal{G}_{\eps}
      \left( \m^{\eps} \right) + \frac{\nu}{2} \Ddelta \right)
    \geqslant \tilde{\mathcal{G}}_0 \left( \m \right) . \end{align}
  {\tmem{}}{\noindent}(iii) ({\tmem{$\Gamma \textrm{\text{-}limsup}$}}
  inequality) We proceed as in the proof of
  Theorem~\ref{thm:thmKS}. Let $\m \in H_0 ( \Om; \Stwo^2 )$ be such
  that $\tilde{\mathcal{G}}_0 \left( \m \right) < + \infty$ and,
  without loss of generality, $\m =\tmmathbf{e}_3$ on $\partial
  \Om$. We take $\bar{\eps} > 0$ and extend $\m$ to
  $\widetilde{\m} \in H^1 \left( \Om_{\bar{\eps}}, \Stwo^2 \right)$ by
  setting $\widetilde{\m} = \tmmathbf{e}_3$ in
  $\Om_{\dbar} \backslash \Om$. For every $\eps < \bar{\eps}$ we now
  define $\m^{\eps} = \widetilde{\m}$ in $\Om_\eps$ and $\m^\eps = 0$
  outside $\Om_\eps$. It is clear that $\m^{\eps} \in \myspace$ and
  satisfies $\m^\eps \rightarrow \m$ strongly in $L^2( \RR^2; \RR^3)$
  as $\eps \rightarrow 0$. Moreover, due to the fact that
  $\m^{\eps} = \m$ in $H^1 \left( \Om; \Stwo^2 \right)$, we have
  $\m^{\eps} =\tmmathbf{e}_3$ on $\partial \Om$. Noting that in this
  case the inequality in \eqref{eq:ineqgammaliming} is actually an
  equality, and using the fact that
  \begin{align} \int_{\Om_{\eps}} \eta_{\eps }^2 \left| \nabla \m^{\eps} \right|^2 
     \dd \rr = \int_{\Od} \eta_{\eps }^2  \left| \nabla \widetilde{\m}
     \right|^2 \dd \rr + \int_{\Om} \left| \nabla \m \right|^2  \dd \rr \; =
     \; \int_{\Om} \left| \nabla \m \right|^2  \dd \rr, \end{align}
  we obtain
  \begin{align} \limsup_{\eps \rightarrow 0}  \left( \mathcal{G}_{\eps} \left(
     \m^{\eps} \right) + \frac{\nu}{2} \Ddelta  \right) =
     \tilde{\mathcal{G}}_0 \left( \m  \right) . \end{align}
  This completes the proof.
\end{proof}

\begin{remark}
  An examination of the proof of Theorem \ref{thm:nonlocengy} shows
  that
  \begin{align}
    | \ln
    \eps | \int_{\partial \Om}  \left( \msq \cdot \nnn \right)^2  \dd
    \mathcal H^1(\sigma)  \to 0, \qquad | \ln \eps |
    \int_{\partial \Om} \left( 1 - | \mps |^2 \right) \dd \mathcal
    H^1(\sigma) \to 0, 
  \end{align}
  as $\eps \to 0$ for any sequence of minimizers $\m^\eps$ of
  $\mathcal G_\eps$.
\end{remark}



\bibliographystyle{acm}
\bibliography{mura}

\end{document}